\documentclass{article}

\usepackage{varioref}   
\usepackage[all]{xy}    
\usepackage{graphics}
\usepackage{graphicx}
\usepackage[centertags]{amsmath}
\usepackage{amssymb}
\usepackage{amsthm}

\newtheorem{thm}{Theorem} 

\newtheorem{prop}{Proposition} 
 
\newtheorem{cor}[prop]{Corollary}

\newtheorem{rem}[prop]{Remark}

\newcounter{diagramm}   

\newcommand{\captionfigure}{\refstepcounter{diagramm}{\it Figure }\arabic{diagramm}}

\newcommand{\CC}{\mathbb{C}}

\newcommand{\RR}{\mathbb{R}}
\newcommand{\ZZ}{\mathbb{Z}}

\newcommand{\HH}{\mathbb{H}} 
 
\newcommand{\AAA}{\mathbb{A}}
\renewcommand{\AA}{\mathbb{A}}
\newcommand{\PP}{\mathbb{P}}

\newcommand{\GT}{\mbox{GT}}
\newcommand{\GTdach}{\widehat{\GT}}

\newcommand{\aut}{\mbox{Aut}}

\newcommand{\torus}{E^{\star}}

\newcommand{\affplus}{\mbox{Aff}^+}

\newcommand{\slzwei}{\mbox{SL}_2}

\newcommand{\abar}{\bar{c}}
\newcommand{\drehung}{c}
\newcommand{\trn}{E^{[n]}}
\newcommand{\tr}{E^{[1]}}

\newcommand{\pgldrei}{\mbox{PGL}_3}
\newcommand{\equ}{E^{[1]}}

\newcommand{\PPP}{\hspace*{-2mm}\begin{array}{l}
                   \vspace*{.6mm}
                   \dddot{\PP}
                   \end{array} \hspace*{-2mm}
                 }

\begin{document}



\title{An extraordinary origami curve}

\date{}
\author{{\it Frank Herrlich
 \footnote{
 e-mail: {\sf herrlich@math.uni-karlsruhe.de}, 
    },
    Gabriela Schmith\"usen}
  \footnote{
   e-mail: {\sf schmithuesen@math.uni-karlsruhe.de}}\\[3mm]
     \small
     Mathematisches Institut II, 
       Universit\"at Karlsruhe, 76128 Karlsruhe, Germany}

\maketitle
\begin{abstract}
\noindent
We study the origami $W$ defined by the quaternion group of order 8 and its Teichm\"uller
curve $C(W)$ in the moduli space $M_3$. We prove that $W$ has Veech group SL$_2(\mathbb Z)$, 
determine the equation of the family over $C(W)$ and find several further properties. 
As main result we obtain infinitely many origami curves in $M_3$ that intersect $C(W)$. 
We present a combinatorial description of these origamis.
\end{abstract}

\vskip 8mm

Origami curves are certain special Teichm\"uller curves in some moduli space of
curves. They are obtained from an unramified covering of a once punctured
torus, see Section \ref{defs} for a precise definition.\\[2mm] 
First examples of Teichm\"uller curves were already given by Veech in 
\cite{V}. In recent years they have attracted a lot of attention, partly because 
of their relation to rational billiards, see e.g.~\cite{mcm} and references therein; 
a survey of  examples defining primitive Teichm\"uller curves can be found 
in \cite{HuSc}. Another topic of interest is the action of the absolute Galois group 
of the rationals on them as discussed in \cite{lo} and \cite{M2}.\\[2mm] 
In some respects Teichm\"uller curves arising via origamis are more accessible 
than general ones. Precisely for them, the Veech group (as defined
in Section \ref{defs}) is a subgroup of $\slzwei(\ZZ)$, cf.~\cite{gj}. It 
can be determined for each origami  explicitly, cf. \cite{sm}. 
Nevertheless, it is still difficult to approach the question
how their 
Teichm\"uller curves are located in the moduli space. There are  
only a few origami curves for which explicit equations
have been found so far, cf. \cite{H} and \cite{M2}.\\[2mm]
In this note we present an extraordinary origami curve in genus 3. It
is the smallest nontrivial example of a normal origami having as Veech group
the full group $\slzwei(\ZZ)$, see Prop.~\ref{sl2z}.
The Jacobian of the
associated family of curves has a two-dimensional fixed part, see
Proposition~\ref{jac}. M.~M\"oller has observed that this implies that our origami
curve is also a Shimura curve; recently he proved that it is the only algebraic
curve in a moduli space of curves of genus at least 2
which is at the same time a Teichm\"uller curve and a
Shimura curve, cf.\ \cite{st}. An explicit equation for the associated family 
of curves is given in Proposition~\ref{fam}. This family of curves has been studied by several authors, and
some of the results of Section~\ref{wms} have been known previously,
cf.~\cite{gi}, \cite{gu}, \cite{kk};
nevertheless the use of the ``origami'' structure allows for new proofs 
and makes the exposition considerably more elementary.\\[2mm]
A very remarkable property of this origami curve, and the main new result of
this paper, is the fact that it intersects infinitely many other origami
curves, see Theorem \ref{thethm} in Section~\ref{sectiondrei}. To our knowledge this 
is the first example of 
origami curves that intersect in moduli space.\\[2mm] 
The decomposition of the
Jacobian gives a second map onto an elliptic curve; if this map
is ramified only over torsion points it can be made into an origami. In Section~\ref{main}
we develop explicit formulas for these maps to determine when this particular
type of ramification occurs.\\[2mm]
We give a combinatorial description of the infinitely many origamis that intersect our origami
curve, see Proposition \ref{schnittorigamis} in Section~\ref{dq}: They can be obtained by suitably glueing two copies
of an $n\!\times\!n\,$-$\,$square (for different~$n$).\\[5mm]
{\bf Acknowledgement}:
   The results in the first section were obtained together
   with Martin M\"oller. We would like to thank him for allowing us to include
   this part. When he visited Karlsruhe and we discovered one remarkable property
   of this origami after the other, we called it ``eierlegende Wollmilchsau''. We
   still call it $W$ in this paper.

\section{The origami $W$}
\label{wms}

We start in \ref{defs} with a brief introduction to origamis. 
In \ref{quat} we introduce the origami $W$, to which the title of this 
paper refers, and prove that it has Veech group $\slzwei(\ZZ)$. In \ref{G} --  \ref{cuspsection} 
we examine properties of $W$: We describe its automorphism group 
(\ref{G}), give algebraic equations for the corresponding family of curves 
(\ref{cw}), find the Jacobian (\ref{jacobi}) and determine the only cusp of its
origami curve (\ref{cuspsection}).
 
\subsection{Origamis and their Veech groups}
\label{defs}
Here we state the basic definitions and facts that we will use in this article.
We follow mostly the notations for origamis used in \cite{sm}. 
A more detailed introduction to translation surfaces and Veech groups can be 
found e.g.~in \cite{gj} and \cite{HuSc}. Some interesting facts about origamis
are contained e.g.~in \cite{le} and \cite{M2}. For an introduction to Teichm\"uller
geodesics and Teichm\"uller curves we refer to \cite{eg} and \cite{mcm}.\\[2mm] 
An {\it origami} is a finite covering $p:X\to E$ from a closed surface $X$ to
a torus $E$ which is unramified over $\torus:=E-\{0\}$, where $0\in E$ is a
point. Any complex uniformization $E_\tau = \CC/(\ZZ+\tau\ZZ)$ of $E$ defines a
translation structure on $E$ which lifts to a translation structure
on $X^\star:=p^{-1}(\torus)=X- p^{-1}(0)$ (therefore an origami is also called
a {\it square tiled surface}). The translation structure is also obtained from
the holomorphic 1-form $\omega_\tau:=p^*\omega_{E_\tau}$, where
$\omega_{E_\tau}$ is the invariant differential on the elliptic curve
$E_\tau$. If $X_\tau$ denotes the complex structure on $X$ determined by this
translation structure, the map $\tau\mapsto X_\tau$ is an isometric embedding
of the upper half plane into the Teichm\"uller space $T_{g,n}$ (where $g$ is
the genus of $X$ and $n= |p^{-1}(0)|$); such embeddings are called {\it
  Teichm\"uller geodesics} or {\it Teichm\"uller disks}.\\[1mm]
Conversely, any finite collection of euclidean squares together with a pairing
that glues every left edge to a right one and every top edge to a bottom one,
uniquely determines an origami. It is this combinatorial construction that
inspired P.~Lochak in \cite{lo} to use the name ``origami''. In the sequel we shall
describe origamis usually in this way.\\[1mm]
For an origami, and more generally for a translation surface $X^\star$ of finite 
volume,
the {\it Veech group} $\Gamma(X^\star)\subset \slzwei(\RR)$ is defined as follows: Let
$\affplus(X^\star)$ be the group of orientation preserving diffeomorphisms of $X^\star$
which are affine with respect to the given translation structure. The linear part of
$\varphi\in\affplus(X^\star)$ is the same on every chart and gives rise to a group
homomorphism der$:\affplus(X^\star)\to \slzwei(\RR)$; then $\Gamma(X^\star)$ is defined as
its image. It is a discrete subgroup of $\slzwei(\RR)$, see \cite{V}.\\[1mm]
The translation surfaces corresponding to the points on a Teichm\"uller disk are related 
by affine diffeomorphisms. Their Veech groups are conjugated in $\slzwei(\RR)$ and it 
is therefore sufficient to consider the Veech group of a base point. 
In the special case of an origami $O$ we choose $\tau = i$ and denote  
$\Gamma(O) := \Gamma(O_i)$. It is known that this is a subgroup of $\slzwei(\ZZ)$ of 
finite index 
(see \cite{sm} for a proof and also for
more details about origamis and Veech groups). Therefore $\Gamma(O)$ acts on the
upper half plane $\HH$ as a lattice, and the orbit space $\HH/\Gamma(O)$ is a
nonsingular affine curve. It is the normalization of the image $C(O)$ of the
Teichm\"uller geodesic defined by $O$ in the moduli space $M_{g,n}$. The
algebraic curve $C(O)$ is called the {\it origami curve} associated with
$O$. This construction is a special case of the more general concept of a Teichm\"uller curve.

\subsection{The quaternion origami}
\label{quat}
The origami $W$, which is the central object of this paper, can be described as
follows: Let $Q:=\{\pm 1,\pm i,\pm j,\pm k\}$ be the quaternion group. Recall
that the (defining) relations are $i^2=j^2=k^2=-1$, $ij=-ji=k$, $(-1)^2 = 1$.
Thus $ik=-ki=-j$ and $jk=-kj=i$. 
Take eight squares labeled with the elements of $Q$ and glue them
in horizontal (resp.\ vertical) direction such that the right neighbour of the
square labeled $g$ has label $g\cdot i$ and its top neighbour has label
$g\cdot j$. The resulting origami $W$ can be represented as
\begin{center}
\setlength{\unitlength}{.8cm}
\begin{minipage}{13cm}
\begin{picture}(5,3)
\put(0,.5){\framebox(1,1){1}}
\put(1,.5){\framebox(1,1){$i$}}
\put(2,.5){\framebox(1,1){$-1$}}
\put(3,.5){\framebox(1,1){$-i$}}
\put(3,1.5){\framebox(1,1){$-k$}}
\put(4,1.5){\framebox(1,1){$-j$}}
\put(5,1.5){\framebox(1,1){$k$}}
\put(6,1.5){\framebox(1,1){$j$}}

\put(8.8,1.5){or}
\put(11,0){\framebox(1,1){$-j$}}
\put(13,0){\framebox(1,1){$j$}}
\put(11,1){\framebox(1,1){1}}
\put(12,1){\framebox(1,1){$i$}}
\put(13,1){\framebox(1,1){$-1$}}
\put(14,1){\framebox(1,1){$-i$}}
\put(12,2){\framebox(1,1){$k$}}
\put(14,2){\framebox(1,1){$-k$}}

\put(.4,1.4){\scriptsize{/}}
\put(1.4,1.4){\scriptsize{$\backslash$}}
\put(2.2,1.4){\scriptsize{///}}
\put(4.2,1.4){\scriptsize{///}}
\put(5.4,1.4){\scriptsize{$\backslash$}}
\put(6.4,1.4){\scriptsize{/}}
\put(.1,.4){\scriptsize{////}}
\put(1.2,.4){\scriptsize{$\backslash\backslash\backslash$}}
\put(2.3,.4){\scriptsize{//}}
\put(3.3,.4){\scriptsize{$\backslash\backslash$}}
\put(3.2,2.4){\scriptsize{$\backslash\backslash\backslash$}}
\put(4.1,2.4){\scriptsize{////}}
\put(5.3,2.4){\scriptsize{$\backslash\backslash$}}
\put(6.3,2.4){\scriptsize{//}}

\put(11.3,-.1){\scriptsize{//}}
\put(13.4,-.1){\scriptsize{/}}
\put(12.3,.9){\scriptsize{$\backslash\backslash$}}
\put(14.4,.9){\scriptsize{$\backslash$}}
\put(11.4,1.9){\scriptsize{/}}
\put(13.3,1.9){\scriptsize{//}}
\put(12.4,2.9){\scriptsize{$\backslash$}}
\put(14.3,2.9){\scriptsize{$\backslash\backslash$}}

\put(10.83,.28){=}
\put(10.83,.45){=}
\put(11.9,.4){--}
\put(12.83,.4){=}
\put(13.83,.45){=}
\put(13.9,.33){--}
\put(11.9,2.4){--}
\put(12.83,2.4){=}
\put(13.83,2.45){=}
\put(13.9,2.33){--}
\put(14.83,2.28){=}
\put(14.83,2.45){=}

\end{picture}
\begin{center}
\captionfigure \label{orig}\\
Edges are glued if they have the same label or are opposite and unlabeled.
\end{center}
\end{minipage}
\end{center}
\vskip 5mm
The commutator of the generators is $ij(-i)(-j)=-1$, thus of order 2. This
implies that every vertex belongs to 8 squares. Therefore the number of
vertices on $W$ is $(8\cdot 4):8 = 4$. The Euler formula gives
$2-2g=8-16+4=-4$, so the genus of $W$ is 3.
\begin{rem}
The origami map $p:W\to E$, which maps each of the eight squares to the
standard torus $E$, is a normal covering with Galois group $Q$. The elements of $Q$ act as translation
automorphisms on $W$ by left(!) multiplication on the labels of the
squares.
\end{rem}
This follows
from the fact that Figure \ref{orig} can also be seen as the Cayley graph of
$Q$ with respect to the generators $i$ and $j$ (taking the squares as vertices
and the glueings as edges).\\[1mm]
Let $p:W^\star\to E^\star$ be the unramified covering induced by $W$ (i.\,e.\
$W^\star = W-p^{-1}(0)$). By the theorem of the universal covering it defines 
an embedding of $H:=\pi_1(W^\star)$ into $\pi_1(E^\star)\cong F_2$, the free
group on two generators. We choose as generators the horizontal and vertical
simple closed curves $x$ and $y$ on $E^\star$. If $\alpha:F_2\to Q$ denotes
the homomorphism that maps $x$ to $i$ and $y$ to $j$, then $H$ is the kernel of $\alpha$. 
\begin{prop}
\label{sl2z}
(i) $H$ is a characteristic subgroup of $F_2$.\\
(ii) The Veech group $\Gamma(W)$ of $W$ is {\rm SL}$_2(\ZZ)$.
\end{prop}
\begin{proof}
The second statement follows from the first by \cite[Prop.~1]{sm}.\\
To show that $H={\rm ker}(\alpha)$ is preserved by all automorphisms of $F_2$
we shall show that every surjective homomorphism $\beta:F_2\to Q$ has kernel
equal to $H$. For this it suffices to observe that $\beta$ is obtained from
$\alpha$ by composition with an automorphism of $Q$: Let $a:=\beta(x)$ and
$b:=\beta(y)$; since $\beta$ is surjective, $a$ and $b$ are elements of order 4 
satisfying $a\not=b\not=-a$. Then the map $i\mapsto a$, $j \mapsto b$ respects
all the defining relations of $Q$, hence is an automorphism $\sigma$ of $Q$;
by construction, $\beta=\sigma\circ\alpha$.\end{proof}
For a more systematic approach to finding such {\it characteristic origamis} cf. \cite{H}.
\begin{prop}
\label{wmodpm1}
$W/\{\pm 1\}$ is isomorphic to $E$ and the origami map $p:W\rightarrow E$ is the 
composition of the quotient map and multiplication by 2.
\end{prop}
\begin{proof}
Note that $-1=i^2=j^2$ acts on the squares of the origami by translation by 2 in
both horizontal and vertical direction. In particular, $-1$ fixes the four
punctures of $W$, as can be read off from the picture (note that the four
vertices of each square all belong to different punctures on $W$). By the
Riemann-Hurwitz formula this implies that $W/\{\pm 1\}$ is an elliptic
  curve. \\
More precisely, $W/\{\pm 1\}$ is the origami\\[3mm]
\setlength{\unitlength}{.8cm}
\hspace*{2cm}
\begin{minipage}{5cm}
\begin{picture}(5,3.5)
\put(3,0){\framebox(1,1){$-j$}}
\put(3,1){\framebox(1,1){1}}
\put(4,1){\framebox(1,1){$i$}}
\put(4,2){\framebox(1,1){$k$}}

\put(6.5,1.6){which is}
\put(6.5,1){equivalent to}

\put(10,.5){\framebox(1,1)}
\put(11,.5){\framebox(1,1)}
\put(10,1.5){\framebox(1,1)}
\put(11,1.5){\framebox(1,1)}

\put(2.9,.4){=}
\put(3.9,.4){--}
\put(3.9,2.4){--}
\put(4.9,2.4){=}
\put(3.3,-.1){\scriptsize{//}}
\put(3.3,1.9){\scriptsize{//}}
\put(4.4,.9){\scriptsize{/}}
\put(4.4,2.9){\scriptsize{/}}

\put(9.9,1.9){=}
\put(9.9,.9){--}
\put(11.9,.9){--}
\put(11.9,1.9){=}
\put(11.4,.4){\scriptsize{//}}
\put(11.4,2.4){\scriptsize{//}}
\put(10.5,.4){\scriptsize{/}}
\put(10.5,2.4){\scriptsize{/}}

\end{picture}
\end{minipage}
\begin{center}
\captionfigure \label{halb}
\end{center}
This shows that the origami map $p:W\to E$ is the composition of the quotient
map $W\to W/\{\pm 1\}$ and multiplication by 2.\\ 
The second step also follows from the fact that $Q/\{\pm 1\}$ is isomorphic to
the Klein four group $V_4$.
\end{proof} 

\subsection{The automorphism group of $W$}
\label{G}
We shall examine the group $G := \aut(W)$ of affine holomorphic automorphisms living on all 
translation structures defined by $W$. Any such automorphism descends to an affine map on the elliptic 
curve with derivative $I$ or $-I$, where $I$ is the identity matrix (see e.g. \cite{sm}). Obviously,  
$Q \subseteq G$ contains the automorphisms with derivative $I$ (compare \ref{quat}). 
In this section we give the group structure of the whole group $\aut(W)$.\\[1mm] 
$W$ admits a further automorphism $\sigma$ that is not in $Q$: the canonical involution 
(i.e.\ multiplication by $-1$) on $E$ can be lifted to $W$. This follows from the fact that
$-I$ is in the Veech group of $W$. Explicitly, we describe $\sigma$ by rotating each square
around its center, and then glueing the same edges as before:
\\[3mm]
\hspace*{2cm}
\setlength{\unitlength}{.75cm}
\begin{minipage}{5cm}
\begin{picture}(5,3)

\put(0,.5){\framebox(1,1){1}}
\put(1,.5){\framebox(1,1){$i$}}
\put(2,.5){\framebox(1,1){$-1$}}
\put(3,.5){\framebox(1,1){$-i$}}
\put(3,1.5){\framebox(1,1){$-k$}}
\put(4,1.5){\framebox(1,1){$-j$}}
\put(5,1.5){\framebox(1,1){$k$}}
\put(6,1.5){\framebox(1,1){$j$}}

\put(.4,1.4){\scriptsize{/}}
\put(1.4,1.4){\scriptsize{$\backslash$}}
\put(2.2,1.4){\scriptsize{///}}
\put(4.25,1.4){\scriptsize{///}}
\put(5.4,1.4){\scriptsize{$\backslash$}}
\put(6.4,1.4){\scriptsize{/}}
\put(.15,.4){\scriptsize{////}}
\put(1.2,.4){\scriptsize{$\backslash\backslash\backslash$}}
\put(2.3,.4){\scriptsize{//}}
\put(3.3,.4){\scriptsize{$\backslash\backslash$}}
\put(3.2,2.4){\scriptsize{$\backslash\backslash\backslash$}}
\put(4.2,2.4){\scriptsize{////}}
\put(5.3,2.4){\scriptsize{$\backslash\backslash$}}
\put(6.3,2.4){\scriptsize{//}}
\put(3.2,1.4){\scriptsize{$\backslash\backslash\backslash\backslash$}}
\put(-.17,.73){=}
\put(-.17,.9){=}
\put(.9,.9){--}
\put(1.83,.9){=}
\put(2.9,.9){--}
\put(2.9,.81){--}
\put(2.9,.73){--}
\put(3.83,.9){=}
\put(3.83,.73){=}

\put(2.83,1.73){=}
\put(2.83,1.9){=}
\put(3.9,1.9){--}
\put(4.83,1.9){=}
\put(5.9,1.9){--}
\put(5.9,1.81){--}
\put(5.9,1.73){--}
\put(6.83,1.9){=}
\put(6.83,1.73){=}

\put(7,.8){\it rotate}
\put(7.2,.5){$\longrightarrow$}

\put(9,0){\framebox(1,1){1}}
\put(10.5,0){\framebox(1,1){$i$}}
\put(12,0){\framebox(1,1){$-1$}}
\put(13.5,0){\framebox(1,1){$-i$}}
\put(10.5,1.5){\framebox(1,1){$-k$}}
\put(12,1.5){\framebox(1,1){$-j$}}
\put(13.5,1.5){\framebox(1,1){$k$}}
\put(15,1.5){\framebox(1,1){$j$}}

\put(8.9,.4){--}
\put(9.83,.45){=}
\put(9.83,.28){=}
\put(10.33,.4){=}
\put(11.4,.4){--}
\put(11.9,.33){--}
\put(11.83,.45){=}
\put(12.83,.4){=}
\put(13.33,.28){=}
\put(13.33,.45){=}
\put(14.33,.45){=}
\put(14.4,.33){--}
\put(10.4,1.9){--}
\put(11.33,1.78){=}
\put(11.33,1.95){=}
\put(11.83,1.9){=}
\put(12.9,1.9){--}
\put(13.4,1.83){--}
\put(13.33,1.95){=}
\put(14.33,1.9){=}
\put(14.83,1.78){=}
\put(14.83,1.95){=}
\put(15.9,1.83){--}
\put(15.83,1.95){=}
\put(9.4,-.1){\scriptsize{/}}
\put(10.9,-.1){\scriptsize{$\backslash$}}
\put(12.25,-.1){\scriptsize{///}}
\put(13.75,-.1){\scriptsize{$\backslash\backslash\backslash\backslash$}}
\put(9.25,.9){\scriptsize{////}}
\put(10.75,.9){\scriptsize{$\backslash\backslash\backslash$}}
\put(12.3,.9){\scriptsize{//}}
\put(13.8,.9){\scriptsize{$\backslash\backslash$}}
\put(10.75,1.4){\scriptsize{$\backslash\backslash\backslash$}}
\put(12.2,1.4){\scriptsize{////}}
\put(13.8,1.4){\scriptsize{$\backslash\backslash$}}
\put(15.3,1.4){\scriptsize{//}}
\put(10.7,2.4){\scriptsize{$\backslash\backslash\backslash\backslash$}}
\put(12.25,2.4){\scriptsize{///}}
\put(13.9,2.4){\scriptsize{$\backslash$}}
\put(15.4,2.4){\scriptsize{/}}

\end{picture}
\end{minipage}
\vskip 5mm
\hspace{4cm}
\begin{minipage}{5cm}
\begin{picture}(5,3)

\put(-2,1.8){\it glue}
\put(-2,1.5){$\longrightarrow$}
\put(0,.5){\framebox(1,1){1}}
\put(1,.5){\framebox(1,1){$-i$}}
\put(2,.5){\framebox(1,1){$-1$}}
\put(3,.5){\framebox(1,1){$i$}}
\put(3,1.5){\framebox(1,1){$-k$}}
\put(4,1.5){\framebox(1,1){$j$}}
\put(5,1.5){\framebox(1,1){$k$}}
\put(6,1.5){\framebox(1,1){$-j$}}

\put(.4,.4){\scriptsize{/}}
\put(1.2,.4){\scriptsize{$\backslash\backslash\backslash\backslash$}}
\put(3.4,.4){\scriptsize{$\backslash$}}
\put(2.2,.4){\scriptsize{///}}
\put(6.25,2.4){\scriptsize{///}}
\put(5.4,2.4){\scriptsize{$\backslash$}}
\put(4.4,2.4){\scriptsize{/}}
\put(.2,1.4){\scriptsize{////}}
\put(3.2,1.4){\scriptsize{$\backslash\backslash\backslash$}}
\put(2.3,1.4){\scriptsize{//}}
\put(1.3,1.4){\scriptsize{$\backslash\backslash$}}
\put(3.2,1.4){\scriptsize{$\backslash\backslash\backslash$}}
\put(6.2,1.4){\scriptsize{////}}
\put(5.3,1.4){\scriptsize{$\backslash\backslash$}}
\put(4.3,1.4){\scriptsize{//}}
\put(3.2,2.4){\scriptsize{$\backslash\backslash\backslash\backslash$}}
\put(.83,.73){=}
\put(.83,.9){=}
\put(-.1,.9){--}
\put(2.83,.9){=}
\put(1.9,.9){--}
\put(1.9,.81){--}
\put(1.9,.73){--}
\put(3.9,.9){--}

\put(3.83,1.73){=}
\put(3.83,1.9){=}
\put(2.9,1.9){--}
\put(5.83,1.9){=}
\put(4.9,1.9){--}
\put(4.9,1.81){--}
\put(4.9,1.73){--}
\put(6.9,1.9){--}

\end{picture}
\end{minipage}
\begin{center}
\captionfigure \label{dreh}
\end{center}
The fixed points of $\sigma$ are the centers of the squares with labels
$1,-1,k,$ and $-k$. So $\sigma$ can also be described as rotation of $W$
around the center of (one of) these squares. \\
Again by Riemann-Hurwitz we find that $W/\!\!<\!\!\sigma\!\!>$ has genus 1. We shall
determine this elliptic curve below in Section \ref{jacobi}.\\[3mm]
From the picture of the
origami we see that $\tau:=i\sigma$ has order 2, and that the same holds for
$\rho:=j\sigma$. This implies $i\sigma(-i) = -\sigma$ and $j\sigma(-j) =
-\sigma$. As a consequence, $\sigma$ commutes with $k$. The element
$\drehung:=k\sigma$ therefore satisfies
$$\drehung^2 = (k\sigma)^2 = k\sigma\sigma k = k^2 = -1.$$
In particular, $\drehung$ is of order 4; it commutes with $k$ and $\sigma$, and also
with $i$ and $j$, thus it is in the center of $G$. With this information at hand
we have determined the group structure of $G$:
\begin{prop}\label{fix}
\mbox{\bf a)} The group $G = \aut(W)$ is generated by $Q$ and $\sigma$ and has order 16. It consists of eight elements of order $4$,
namely $\pm i,\pm j, \pm k$ and $\pm \drehung$, and seven elements of order $2$, namely
$-1$ and $\pm \sigma,\pm\tau$ and $\pm\rho$. In the list of Hall and Senior \cite{HS} it belongs to the
family $\Gamma_2$ and is number 8 among the groups of order $16$.\\[1mm]
\mbox{\bf b)} 
The automorphisms in $Q$ descend via $p$ to the identity on $E$, whereas 
those in $G\backslash Q$ descend to involutions. 
On the degree 2 quotient $W/\{\pm 1\} \cong E$ (see Proposition~\ref{wmodpm1}) 
the elements of $Q$ induce translations and those in $G \backslash Q$ 
induce involutions.\\[1mm]
\mbox{\bf c)}
The elements of $G$ have the following fixed points in $W$:
\begin{itemize}
\item[] $\pm i, \pm j$ and $\pm k$ have no fixed points.
\vspace{-2mm}
\item[] $\sigma$ fixes the centers of the squares 1, $-1,k,-k$.
\vspace{-2mm}
\item[] $-\sigma$ fixes the centers of the squares $i,-i,j,-j$.
\vspace{-2mm}
\item[] $\tau$ fixes the centers of the vertical edges between 1 and $i$, $-1$
  and $-i$, $-j$ and $k$, $j$ and $-k$.
\vspace{-2mm}
\item[] $-\tau$ fixes the centers of the remaining vertical edges.
\vspace{-2mm}
\item[] $\rho$ fixes the centers of the horizontal edges between 1 and $j$,
  $i$ and $k$, $-1$ and $-j$,  $-i$ and $-k$.
\vspace{-2mm}
\item[] $-\rho$  fixes the centers of the remaining horizontal edges.
\vspace{-2mm}
\item[] $\drehung$ (and hence also $-1$) fixes the vertices of the squares.
\end{itemize}
\end{prop}
\subsection{The origami curve $C(W)$}
\label{cw}
By Proposition~\ref{fix}, $C(W)$ is a 1-parameter
family of curves of genus 3 with automorphism group $G$. This uniquely
determines this family, see \cite{kk}. Using the origami structure we can
give an independent direct proof:
\begin{prop}
\label{fam}
\mbox{\bf a)} The origami curve $C(W)$ is the image in $M_3$ of the 1-parameter family of
smooth plane curves $W_{\lambda}$ with affine equation
$$y^4=x(x-1)(x-\lambda) \quad \mbox{ with } \lambda \in \PPP := \PP^1(\CC) - \{0,1,\infty\}.$$
The quotient map $W_{\lambda} \rightarrow W_{\lambda}/\{\pm 1\}=E_{\lambda}: y^2 = x(x-1)(x- \lambda)$ 
is given by 
$$(x,y) \mapsto (x,y^2).$$
\mbox{\bf b)} $C(W)$ is nonsingular and isomorphic to $\HH/\Gamma(W)=\AA^1$.
\end{prop}
\begin{proof}
a) From the Riemann-Hurwitz formula and the fact that $\drehung$ has four fixed
points on $W$ we conclude $4 = 4\cdot(2g(W/\!\!<\!\drehung\!>\!)-2)+4\cdot3$
which implies that $W/\!\!<\!\drehung\!>$ is a rational curve. Thus $W$ is a cyclic covering of order 4
of the projective line, totally ramified over four distinct
points. Normalizing three of them to $0,1$ and $\infty$ we see that $W$ has
the equation
$$y^4 = x^{\varepsilon_0}(x-1)^{\varepsilon_1}(x-\lambda)^{\varepsilon_\lambda},$$
where the parameter $\lambda$ is different from 0, 1 and $\infty$ and depends on (or rather determines) the complex
structure of $W$, and the exponents $\varepsilon_i$ can be 1 or 3 (2 is
excluded because the covering is totally ramified over all four points). Replacing $y$
by $y^{-1}x(x-1)(x-\lambda)$ if $\varepsilon_0=3$ we may assume
$\varepsilon_0=1$. The values of $\varepsilon_1$ and $\varepsilon_\lambda$ are
determined by the monodromy action of $\drehung$ on the loops around $0,1$ and
$\lambda$: from the origami we read off that $\drehung$ acts by counterclockwise
rotation at all four fixed points. This implies that the exponent is the same
at all ramification points and thus proves the first statement.\\
The map given in the claim is compatible with the automorphism $-1$ of $W_{\lambda}$, 
which is   $-1:(x,y) \mapsto (x,-y)$. 
Therefore it is the quotient map.\\[1mm]
b) From Proposition~\ref{sl2z} we know that the Veech group $\Gamma(W)$ is
$\slzwei(\ZZ)$, thus $\HH/\Gamma(W)$ is isomorphic to the affine line. On the
other hand, $C(W)$ is the image of the map $\PPP\to M_3$
that sends $\lambda$ to the isomorphism class of the plane quartic $W_\lambda$. 
Since an embedding of a genus three
curve into the projective plane is necessarily the canonical embedding, two
such curves are isomorphic as abstract curves if and only if they are isomorphic
as embedded curves. This means that for $\lambda$,
$\lambda'\in\PPP$, $W_\lambda$ and $W_{\lambda'}$ define the
same point in $M_3$ if and only if there is a projective change of coordinates
that maps $W_\lambda$ onto $W_{\lambda'}$. As in the case of the Legendre
family $y^2=x(x-1)(x-\lambda)$ of elliptic curves such a projectivity exists
if and only if $\lambda'$ is obtained from $\lambda$ by a M\"obius
transformation that permutes the three points 0, 1, $\infty$. This shows that
$C(W)$ is isomorphic to $\PPP/S_3\cong\AA^1$.
\end{proof}
\subsection{The Jacobian of $W$}
\label{jacobi}
We first determine the quotient of $W$ by $\sigma$: Since $\sigma$ commutes
with $\drehung$ but is not in the center, the normalizer of $\sigma$ in $G$ is a group of order
8, isomorphic to $\ZZ/2\ZZ\times\ZZ/4\ZZ$. This shows that $c$ descends to an automorphism 
$\bar{\drehung}$ of order 4 on the elliptic curve $E_\sigma:=W/\!\!<\!\sigma\!>$. Since $\bar{\drehung}$
has a fixed point (see Proposition \ref{fix}), it is not a translation. The only elliptic curve with such an
automorphism is $E_{-1}$ which has equation $y^2=x^3-x$ and $j$-invariant
1728. By symmetry the same reasoning also holds for the other involutions:
\begin{prop}
\label{wmodsigma}
For any $\lambda \in \PPP$, the quotient of $W_{\lambda}$ by each of the involutions 
$\pm \sigma,\pm\tau$ and
$\pm\rho$ is the elliptic curve $E_{-1}$.
\end{prop}
The following result is contained in
\cite[Thm.\,7.3]{Gu1}, where also the degree of the isogeny is determined. It
is also a special case of a much more general result in \cite{W}. We include
an elementary proof which makes essential use of the origami structure.
\begin{prop}
\label{jac}
For every $\lambda\in \PPP$, the Jacobian of $W_\lambda$ is
isogenous to $E_\lambda\times E_{-1}\times E_{-1}$. 
\end{prop}
\begin{proof}
Fix a quotient morphism $\varphi: W_\lambda\to
W_\lambda/\!\!<\!\sigma\!>\,\cong E_{-1}=:E$. Since $-\sigma
= i\sigma(-i)$ we see that $\psi:=\varphi\circ i:W_\lambda\to E$ is a
quotient map for the action of $-\sigma$. \\[1mm]
Let $\beta:{\rm Jac}\,W_\lambda\to E\times E$ be the homomorphism (of abelian
varieties) induced by $(\varphi,\psi):W_\lambda\to E\times E$. We claim that
$\beta$ is surjective. Together with the quotient map $W_\lambda\to
W_\lambda/\{\pm1\}=E_\lambda$ from Proposition~\ref{wmodpm1} we then obtain the
desired isogeny.\\[1mm] 
To prove the claim let $P_0$ and $P_1$ be the centers of the squares 1 and
$-1$, and $Q_0$, $Q_1$ the centers of the squares $i$ and $-i$. Then
$P_0$ and $P_1$ are fixed by $\sigma$ and $\sigma(Q_0)=Q_1$, while $-\sigma$
fixes 
$Q_0$ and $Q_1$ and interchanges $P_0$ and $P_1$. Moreover $P_0,Q_0,P_1,Q_1$
are (in this order) an orbit under $i$. Denoting $D$ and $D'$ the classes in
Jac$\,W_\lambda$ of the divisors $P_1-P_0$ and $Q_0-Q_1$, respectively, we find
$\beta(D)=(P,0)$ and $\beta(D')=(0,P)$ where
$P=\varphi(P_1)-\varphi(P_0)\not=0$ on $E$. Note that in fact
$\psi(Q_0)-\psi(Q_1)=\psi\circ i(P_0)-\psi\circ
i(P_1)=\varphi(-P_0)-\varphi(-P_1)=\varphi(P_1)-\varphi(P_0)=P$.\\[1mm]
Since the image of $\beta$ must be a connected subgroup, to prove surjectivity
of $\beta$ it is enough to show that there is no onedimensional connected
subgroup  of $E\times E$ containing both $(P,0)$ and $(0,P)$ (for some
$P\not=0$). But all such subgroups of $E\times E$ are of the form
$V(a,b)=\{(x,y)\in E\times E: ax+by=0\}$ for some $a,b\in \ZZ[i]$;
moreover $V(a,b)$ is connected if and only if $a$ and $b$ are relatively
prime. Now the hypothesis $(P,0)\in V(a,b)$ implies $a\cdot P=0$, while
$(0,P)\in V(a,b)$ implies $b\cdot P=0$. Together with gcd$(a,b)=1$ this gives
the contradiction $P=0$.
\end{proof}
\subsection{The cusp of $C(W)$}\label{cuspsection}
Since $C(W)$ is isomorphic to $\AA^1$, its closure in $\overline{M}_3$
contains a unique cusp $W_\infty$. 
\begin{prop}
The stable curve corresponding to $W_\infty$ consists of two irreducible
components, both isomorphic to $E_{-1}$, that intersect transversally in two
points. 
\end{prop}
\begin{proof}
$W_\infty$ can be obtained by
contracting the (dotted) central lines of the horizontal cylinders of the
origami as indicated in Figure \ref{cusp}:
\\[3mm]
\setlength{\unitlength}{.8cm}
\begin{center}
\begin{minipage}{5cm}
\begin{picture}(5,2)
\put(0,.5){\framebox(1,1){- - -}}
\put(1,.5){\framebox(1,1){- - -}}
\put(2,.5){\framebox(1,1){- - -}}
\put(3,.5){\framebox(1,1){- - -}}
\put(3,1.5){\framebox(1,1){- - -}}
\put(4,1.5){\framebox(1,1){- - -}}
\put(5,1.5){\framebox(1,1){- - -}}
\put(6,1.5){\framebox(1,1){- - -}}

\put(.4,1.4){\scriptsize{/}}
\put(1.4,1.4){\scriptsize{$\backslash$}}
\put(2.2,1.4){\scriptsize{///}}
\put(4.2,1.4){\scriptsize{///}}
\put(5.4,1.4){\scriptsize{$\backslash$}}
\put(6.4,1.4){\scriptsize{/}}
\put(.1,.4){\scriptsize{////}}
\put(1.2,.4){\scriptsize{$\backslash\backslash\backslash$}}
\put(2.3,.4){\scriptsize{//}}
\put(3.3,.4){\scriptsize{$\backslash\backslash$}}
\put(3.2,2.4){\scriptsize{$\backslash\backslash\backslash$}}
\put(4.1,2.4){\scriptsize{////}}
\put(5.3,2.4){\scriptsize{$\backslash\backslash$}}
\put(6.3,2.4){\scriptsize{//}}

\end{picture}
\end{minipage}\\
\captionfigure \label{cusp}
\end{center}
Clearly we obtain two irreducible components, that intersect in two points and
thus both are nonsingular of genus~1.\\[1mm] 
From Proposition~\ref{jac} we know that the Jacobian 
of $W_\lambda$ contains a subvariety isogenous to $E_{-1}\times E_{-1}$,
independently of $\lambda$. Therefore the Jacobian of $W_\infty$
must also contain such a subvariety. From this we deduce that the irreducible
components of $W_\infty$ are isogenous to $E_{-1}$. Moreover we find from
the picture that the automorphism $\drehung$ of $W$ has two fixed points on
each component of $W_\infty$ and in particular induces an automorphism on each
component. As in the proof of Proposition~\ref{wmodsigma} we conclude that the
components are isomorphic to $E_{-1}$.
\end{proof}
Any stable curve of genus 3 that consists of two elliptic curves
intersecting in two points is hyperelliptic (i.\,e.\ in the closure of the
hyperelliptic locus in $M_3$), see \cite[ex.~3.162 on p.~185]{HM}. In the case of
$W_\infty$ the hyperelliptic involution $h$ can easily be described explicitly
on the origami: $h$ acts on the ``upper'' horizontal cylinder like $\rho$ and
on the 
``lower'' one like $-\rho$; this means that $h$ fixes the centers of all
horizontal edges and acts by rotation by $\pi$ on each of the squares around
these points. Note that $h$ is not an automorphism of $W$ because the actions
of $\rho$ and $-\rho$ on the two horizontal lines separating the cylinders do
not agree. But on each component of $W_\infty$, $h$ acts as the multiplication
by $-1$, and it interchanges the two intersection points. Thus the quotient of
$W_\infty$ by $h$ consists of two projective lines meeting transversally in
one point. Together with the four punctures of $W$, two in each horizontal
cylinder, this is a stable 4-marked curve of genus 0.\\[3mm]
It is also possible, but a bit harder, to determine $W_\infty$ out
of the explicit formula for the family that we obtained in Proposition \ref{fam}. 
The fibre over 0 (or 1 or $\infty$) in this family is the plane curve with
equation $y^2=x^2(x-1)$, which is irreducible with a tacnode at the origin; in
particular it is not a stable curve. The stable reduction theorem \cite[Prop.~3.47]{HM} guarantees that there is a unique stable curve that, after passing
to a branched covering of the base, fits as special fibre into this family. To
find it we first have to desingularize the total space of the family. This
requires three consecutive blow-ups in the singular point. After that the
special fibre contains a nonreduced rational irreducible component (of
multiplicity 2). Taking a double covering which is ramified precisely over the
reduced components of the special fibre transforms the double line into a
reduced elliptic curve $E$. It turns out that the other components of the
special fibre are a second elliptic curve $E'$ which intersects $E$ in two
points, and two projective lines $L_1$, $L_2$ each of which intersects $E'$ in
one point. So far the special fibre is only semistable, but the total space is
nonsingular. Finally we have to contract $L_1$ and $L_2$ to obtain a stable
special fibre, although this unavoidably leads to two singular points in the
total space. The details of this construction are worked out in~\cite{B}.


\section{Other maps to elliptic curves}
\label{main}
In this section we calculate the action of the elements of $G = \aut(W)$ on $W$
\footnote{Such formulas can also be found in \cite{kk} and \cite{gi}}.
We obtain explicit formulas for the quotient maps of $W$
with respect to some of the automorphisms in $G$. In particular we can
determine their ramification points. This will enable us to show in 
Section \ref{sectiondrei} that the quaternion origami curve $C(W)$ intersects 
infinitely many other origami curves.
 
\subsection{Formulas for the automorphisms}
\label{autos}
In Section \ref{G} we described the automorphisms of $W$ by their action on
the origami. Using the equation of Proposition~\ref{fam} for the curve $W_{\lambda}$  
we translate this into formulas which are summarized in Proposition~\ref{formeln}.\\[1mm]
Recall that the central automorphism $\drehung$ of $G$ has order 4 and descends to 
an involution of 
$E_{\lambda} = W_{\lambda}/\{\pm 1\}$, with the quotient map given in 
Proposition~\ref{fam}. Thus $\drehung$ acts in homogenous coordinates by 
\[\drehung :(X:Y:Z) \mapsto (X:iY:Z).\]
Here we made a choice between $\drehung$ and $-\drehung$ which 
corresponds to the choice between the given map and the map 
$(X:Y:Z) \mapsto (X:-iY:Z)$.\\[2mm]
The four vertices of the origami squares are the four fixed points
of $\drehung$:  
\begin{eqnarray*}
P_0 = (0:0:1), && P_1 = (1:0:1),\\
P_{\infty} = (1:0:0), && P_{\lambda} = (\lambda:0:1).
\end{eqnarray*}
Their images on $E_\lambda$ are the points $Q_0=P_0$, $Q_1=P_1$,
$Q_\infty=(0:1:0)$ and $Q_\lambda=P_\lambda$. Once $Q_\infty$ is fixed as the
origin for the group law on $E_\lambda$, these are the points of order 2, see
Proposition~\ref{wmodpm1}. We choose the translation structure on $E_\lambda$ such
that the horizontal and vertical neighbours of $Q_0$ are $Q_1$ and $Q_\infty$, respectively.
Lifting the structure to $W_{\lambda}$ the vertices of the squares are labeled as follows:\\
\setlength{\unitlength}{1cm}
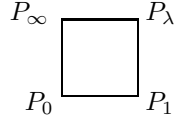
\begin{center}
\begin{minipage}{4cm}
\begin{picture}(2,2)
\put(0.5,1){\framebox(1,1)}
\put(0,.8){$P_0$}
\put(-.2,2){$P_\infty$}
\put(1.6,2){$P_\lambda$}
\put(1.6,.8){$P_1$}
\end{picture}
\end{minipage}
\end{center}
\vspace*{-5mm}
That this choice can be made in a consistent way for the whole family reflects the
fact that the family ${\cal E}:=\{E_\lambda:\lambda\in\PPP\}$
carries a natural (even universal) level 2 structure.\\[2mm]
By our normalization $\sigma$ interchanges $P_0$ with $P_{\lambda}$ and $P_1$
with $P_{\infty}$, see Section~\ref{G}. In the same way, $\rho$ interchanges
$P_0$ with $P_1$ and $P_\infty$ with $P_\lambda$, and $\tau$ interchanges
$P_0$ with $P_\infty$ and $P_1$ with $P_\lambda$.\\[2mm]
Any automorphism in $G$ has to permute these four points, and no
element of $G-<\drehung>$ fixes any of them. Thus $G$ acts through the
normal subgroup $V_4$ of $S_4$.\\
If we represent the automorphisms of $W_{\lambda}$ by matrices in $\pgldrei(\CC)$
with respect to the homogeneous coordinates $X$, $Y$, $Z$, then $\drehung$ is given
by 
$$\tilde\drehung:= \begin{pmatrix}1&0&0\\0&i&0\\0&0&1\end{pmatrix}.$$ 
Since
each $g\in G$ commutes with $\drehung$, it is represented by a matrix $\tilde
g$ such that $\tilde g\tilde \drehung\tilde g^{-1}\tilde
\drehung^{-1}=\mu\cdot I_3$ for some $\mu\in\CC^\times$ and the unit matrix
$I_3$. The determinant of a commutator is 1; this implies $\mu^3=1$. Let $x$
be an eigenvector of $\tilde\drehung$ for the eigenvalue $i$. Then $\tilde
g(x)$ is an eigenvector of $\tilde\drehung$ for the eigenvalue $i\mu^{-1}$,
because  $\tilde\drehung(\tilde g(x))=\mu^{-1}\tilde
g\tilde\drehung(x)=\mu^{-1}\tilde g(ix)$. But the only eigenvalues of
$\tilde\drehung$  are $i$ and 1, hence
$\mu$ must be 1 and so $\tilde g$ commutes with $\tilde\drehung$. Therefore
there is a basis of common eigenvectors of $\tilde\drehung$ and $\tilde
g$. This implies that $\tilde g$  has the form 
\[\phantom{bbbbbbbbbbbbb}\begin{pmatrix}a&0&b\\0&\zeta&0\\c&0&d\end{pmatrix} \qquad
\mbox{for some $\zeta,a,b,c,d \in \CC$.}\]
Now let $\tilde g$ be a matrix representing one of the automorphisms
$\sigma\drehung^\nu$, $\nu=0,1,2,3$, i.\,e.\ $\pm\sigma$ or $\pm k$ (recall that
$c = k\sigma$ and $c^2 = -1$). From the
fact that they all interchange $P_0$ with $P_{\lambda}$ and $P_1$ with
$P_{\infty}$, it is immediate to deduce the conditions $b = -\lambda a$, $c =
a$, $d = -a$. The corresponding automorphism $\alpha_\zeta$ of $\PP^2$ is therefore: 
\begin{eqnarray}
\label{sigma}
\alpha_\zeta (X:Y:Z)= (X-\lambda Z:\zeta Y: X-Z ).
\end{eqnarray}
Finally, $\zeta$ is determined by the fact that $\alpha_\zeta$ has to map $W_\lambda$  onto itself:
\begin{eqnarray*}
\zeta^4 Y^4 =  (X - \lambda Z)(X-Z)Z(1-\lambda)X(1-\lambda) = (1-\lambda)^2XZ(X-Z)(X-\lambda Z) 
\end{eqnarray*}
This gives the condition 
\[\zeta^4 = (1-\lambda)^2.\]
The four possible values of $\zeta$ correspond to the four elements $\pm \sigma$, $\pm k$. To find out which values correspond to the elements of order 2 (i.e. $\pm \sigma$) we look at the square of $\alpha_\zeta$:
\begin{eqnarray}
\label{sigmaquadrat}
\alpha_\zeta^2(X:Y:Z)= (X(1-\lambda): \zeta^2Y : Z(1-\lambda))
\end{eqnarray}
Thus $\alpha_\zeta$ has order 2 if and only if $\zeta^2 = 1 - \lambda$.\\[2mm]
Carrying out the analogous calculations for the other elements of $G$ we
obtain
\begin{prop}
\label{formeln}
For any $\lambda\in\PPP$ the 16 elements of $G$ act on
$W_\lambda$ through the following projective automorphisms
\begin{eqnarray*}
\iota_\nu(X:Y:Z) &=& (X:i^\nu Y:Z),\ \ \nu=0,\dots,3\\
\alpha_\zeta(X:Y:Z) &=& (X-\lambda Z:\zeta Y: X-Z ),\ \ \zeta^4=(1-\lambda)^2\\
\beta_\xi(X:Y:Z) &=& \left(X-Z:\xi Y : \frac{1}{\lambda}X - Z\right),\ \  
\xi^4 = \left(1 - \frac{1}{\lambda}\right)^2\\
\gamma_\eta(X:Y:Z) &=& (\lambda Z:\eta Y : X),\ \ 
\eta^4 = \lambda^2
\end{eqnarray*}
The correspondence with the elements of $G$ is:
\begin{eqnarray*} 
1 = \iota_0,\ \ -1 = \iota_2&\mbox{ and }& \pm c = \iota_\nu\ \ \mbox{for}\ \ \nu = 1,3\\
\pm \sigma = \alpha_\zeta\ \ \mbox{for}\ \ \zeta^2=1-\lambda&\mbox{ and }&\pm
k = \alpha_\zeta\ \ \mbox{for}\ \ \zeta^2=\lambda-1\\
\pm\rho = \beta_\xi\ \ \mbox{for}\ \ \xi^2 = 1-\frac{1}{\lambda}
&\mbox{ and }& 
 \pm i = \beta_xi\ \ \mbox{for}\ \ \xi^2 = \frac{1}{\lambda} - 1\\
 \pm \tau = \gamma_\eta\ \ \mbox{for}\ \ \eta^2 = \lambda
&\mbox{ and }&
 \pm j = \gamma_\eta\ \ \mbox{for}\ \ \eta^2 = -\lambda
\end{eqnarray*}
\end{prop}
\subsection{Quotient by $\sigma$}
\label{proj}
In Proposition~\ref{wmodsigma} we saw that the quotient map with respect to the 
automorphism $\sigma$ is a map of degree 2 from $W_{\lambda}$ to $E_{-1}$.
In this section we obtain an explicit formula for this map $\kappa$ stated 
in Proposition~\ref{kprop}.\\[2mm]
Proposition~\ref{formeln} gives us two possibilities for $\sigma$. We fix a 
square root $\zeta$ of $1-\lambda$
and call $\sigma_\lambda$ the automorphism given by
\[\sigma_\lambda(X:Y:Z) = (X - \lambda Z: \zeta Y: X-Z).\]
The other possible choice is $-\zeta$; it describes the automorphism $(-\sigma)_\lambda$.\\
$\sigma_\lambda$ leaves the affine plane $Y \neq 0$ invariant, and it acts on the affine coordinates $x =\frac{X}{Y}$ and $z = \frac{Z}{Y}$ by
\[ \sigma_\lambda(x,z) = \left(\frac{1}{\zeta}(x - \lambda z),\frac{1}{\zeta}(x-z)\right).\]
The algebra of polynomials in $x$ and $z$ that are invariant under $\sigma_\lambda$ is generated
by 
\begin{eqnarray*}
p'_2 &:=& x + \sigma_\lambda(x) = x\left(1 + \frac{1}{\zeta}\right) - \frac{\lambda}{\zeta}\cdot z\\
\mbox{ and } \quad p'_3 &:=& x\cdot\sigma_\lambda(x) = \frac{1}{\zeta}x\left(x-\lambda z\right).
\end{eqnarray*}
We shall work with 
\begin{eqnarray*}
p_2 &:=& \frac{\zeta}{1+\zeta}\cdot p'_2 =  x - (1 - \zeta)z\\
\mbox{ and } \quad p_3 &:=& \zeta\cdot p'_3 =x(x-\lambda z).
\end{eqnarray*}
The intersection of $W_{\lambda}$ with the affine plane $Y \neq 0$ is the affine curve with the equation 
\[1 = xz(x-z)(x-\lambda z).\]
This polynomial is invariant under $\sigma_\lambda$, hence can be written as a polynomial in $p_2$ and $p_3$: Observe that 
\begin{eqnarray*}
p_2^2 - p_3 &=& x^2 + (1 - \zeta)^2z^2 - 2xz(1-\zeta) - x^2 + \lambda xz\\ 
 &=& (1-\zeta)^2z^2 - xz(2-2\zeta-1+\zeta^2) = -(1-\zeta)^2z(x-z).
\end{eqnarray*}
Thus the affine equation for $W_{\lambda}$ can be rewritten as 
\begin{eqnarray}
\label{stern}
 1 = p_3(p_2^2 - p_3)\cdot(-(1-\zeta)^2)^{-1}.
\end{eqnarray}
Considering $p_2$ and $p_3$ as the affine coordinates for $\AAA^2/\!<\!\sigma_\lambda\!>$, (\ref{stern}) is also the affine equation for $W_{\lambda}/\!<\!\sigma_\lambda\!>$.\\
To bring this cubic equation into Weierstrass form we first pass to the projective closure by homogenizing (\ref{stern}):
\[p_2^2p_3 = -(1-\zeta)^2w^3 + wp_3^2.\]
Choose a square root $\omega$ of $1-\zeta$ and rewrite this equation in the
coordinates $\tilde{p}_2 = \omega p_2$, $\tilde{p}_3 = p_3$, $\tilde{w}:= -(1 - \zeta)w$:
\[ \tilde{p}_2^2\tilde{p}_3 = (1 - \zeta)p_2^2p_3 = -(1-\zeta)^3w^3 + (1-\zeta)wp_3^2
    = \tilde{w}^3 - \tilde{w}\tilde{p}_3^2.\]
This is the standard Weierstrass form of $E_{-1}$.
Putting everything together we get
\begin{prop} \label{kprop}
The quotient map
$\kappa = \kappa_{\sigma_\lambda}: W_{\lambda} \rightarrow W_{\lambda}/\!<\!\sigma_\lambda\!> = E_{-1}$
is given by 
\[\kappa(X:Y:Z) = \left(-(1-\zeta)Y^2:\omega Y(X -
(1-\zeta)Z):X(X-\lambda Z)\right)\]
for $(X:Y:Z)\not\in\{P_0,P_\lambda\}$ and
$\kappa(P_0)=\kappa(P_\lambda)=(0:1:0)$.\\
$\kappa$ is determined by $\lambda$ uniquely up to the sign of $\zeta=\sqrt{1-\lambda}$ (which
determines the choice between $\sigma_\lambda$ and $(-\sigma)_\lambda$), and the sign of
$\omega=\sqrt{1-\zeta}$ (a change of which corresponds to a composition with the
automorphism $[-1]$ on $E_{-1}$).
\end{prop}
\begin{proof}
Considered on the appropriate affine parts of $\PP^2$, the affine
quotient map $(x,y) \mapsto (p_2,p_3)$ is given by $(\frac{X}{Y}:1:\frac{Z}{Y})
\mapsto (1:\frac{p_2}{w}:\frac{p_3}{w})$. The change of coordinates from $w$,
$p_2$, $p_3$ to $\tilde w$, $\tilde p_2$, $\tilde p_3$ transforms this into
$(1:\frac{\omega p_2}{-(1-\zeta)w}:\frac{p_3}{-(1-\zeta)w})$.
Substituting
$\frac{X}{Y}-(1-\zeta)\frac{Z}{Y}$ for $\frac{p_2}{w}$ and
$\frac{X}{Y}(\frac{X}{Y}-\lambda\frac{Z}{Y})$ for $\frac{p_3}{w}$ gives the
formula on $W_\lambda\cap\{Y\not=0\}$. To determine the image of $\kappa$
in the remaining points $P_0$ and $P_\lambda$, replace the third component
$X(X-\lambda Z)$ by $\frac{Y^2}{X(X-Z)}$ and multiply all components
by$\frac{X(X-Z)}{Y}$. Thus on the open set $\{X\not=0,X\not=Z,Y\not=0\}$, we
also have
$\kappa(X:Y:Z)=(-(1-\zeta)YZ(X-Z):\omega(X-(1-\zeta)Z)Z(X-Z):Y^3)$.
This formula gives us $\kappa(P_0)$ and $\kappa(P_\lambda)$.
\end{proof}

\subsection{The other involutions}
\label{others}
To determine the quotient maps $W_\lambda\to E_{-1}$ with respect to the
automorphisms $\rho$ and $\tau$ we proceed as in Section~\ref{proj}.\\[2mm]
Out of the two possibilities for $\rho$ in Proposition~\ref{formeln} we pick one by
fixing a square root $\xi$ of $1-\frac{1}{\lambda}$ and denoting $\rho_\lambda$
the automorphism of $W_\lambda$ which in affine coordinates
$x=\frac{X}{Y}$, $z=\frac{Z}{Y}$ is given by
\[\rho_\lambda(x,z) =
\left(\frac{1}{\xi}\left(x-z\right),\frac{1}{\xi}\left(\frac{1}{\lambda}x-z\right)\right).\]
As in the case of $\sigma$, the other square root $-\xi$ then describes
$(-\rho)_\lambda$. \\[2mm]
This time, the invariants are generated by
\[p_2=x-\frac{1}{\xi+1}z\ \ \mbox{and}\ \ p_3=x(x-z).\]
In these invariants, the right hand side of the equation defining $W_\lambda$
is 
\[xz(x-z)(x-\lambda z)=-p_3(p_2^2-p_3)\cdot(1+\xi)^2\lambda,\]
which leads to the affine equation
\[p_2^2p_3-p_3^2=-\frac{1}{\lambda(1+\xi)^2}\]
for $W_\lambda/\!<\!\rho_\lambda\!>$. As before we homogenize this equation
with a variable $w$ and take a homogeneous change of
coordinates to transform this equation into Weierstrass form:
\[\tilde p_2= \omega_\rho p_2,\ \tilde w=-\omega_\rho^2w,\ \tilde p_3=p_3 \]
with $\omega_\rho^4=\frac{1}{\lambda(1+\xi)^2}$. This gives us the quotient
map $\kappa_\rho=\kappa_{\rho_\lambda}$ stated below in Proposition~\ref{kapparho}.\\[2mm]
Finally $\tau$ acts in affine coordinates by
\[\tau_\lambda(x,z) = \left(\frac{\lambda}{\eta}z,\frac{1}{\eta}x\right),\]
where we fix a square root $\eta$ of $\lambda$ and thus distinguish
$\tau_\lambda$ from $(-\tau)_\lambda$.\\[2mm]
The same reasoning as above gives the quotient map $\kappa_\tau:W_\lambda\to
W_\lambda/\!<\!\tau_\lambda\!> = E_{-1}$ in the following proposition. 
\begin{prop} \label{kapparho}
The quotient maps $\kappa_\rho=\kappa_{\rho_\lambda}:W_\lambda\to
W_\lambda/\!<\!\rho_\lambda\!> = E_{-1}$ and  $\kappa_\tau:W_\lambda\to
W_\lambda/\!<\!\tau_\lambda\!> = E_{-1}$ are given as follows:
\[
\begin{array}{lcl}
\kappa_\rho(X:Y:Z)&=&\left(-\omega_\rho^2Y^2:\omega_\rho 
              Y\left(X-\frac{1}{1+\xi}Z\right):X(X-Z)\right)\\
 \multicolumn{3}{l}{
 \mbox{\quad\qquad for $(X:Y:Z)\not\in\{P_0,P_1\}$ and
 $\kappa_\rho(P_0)=\kappa_\rho(P_1)=(0:1:0)$.}}  \\[2mm]
 \kappa_\tau(X:Y:Z)&=&\left(\frac{i}{1+\eta}Y^2:\omega_\tau Y(X+\eta Z):XZ\right)\\
 \multicolumn{3}{l}{
 \mbox{\quad\qquad
 for $(X:Y:Z)\not\in\{P_0,P_\infty\}$ and $
 \kappa_\tau(P_0)=\kappa_\tau(P_\infty)=(0:1:0)$.}}\\ 
 \multicolumn{3}{l}{
 \mbox{\quad\qquad
  Here $\omega_\tau$ is
  a square root of $(\frac{i}{1+\eta})^3$.}}
\end{array}
\]
$\kappa_\rho$ is determined by $\lambda$ up to the sign of $\xi$ (which gives
the choice between $\rho_\lambda$ and $(-\rho)_\lambda$), and up to
multiplying $\omega_\rho$ by a power of $i$ (which results in composing
$\kappa_\rho$ with a power of the automorphism $\bar\drehung$ of
$E_{-1}$). 
$\kappa_{-\rho}$ is obtained from $\kappa_\rho$ by replacing $\xi$ by $-\xi$.
Similarly $\kappa_{\tau}$ is determined up tp the signs of $\eta$ and $\omega_\tau$. 
\end{prop}
Next we shall show that all these quotient maps are closely related to each
other. First recall that $W_\lambda$ is isomorphic to $W_{\lambda'}$ iff$\lambda'\in \{\lambda,\frac{1}{\lambda},1-\lambda,\frac{1}
{1-\lambda},1-\frac{1}{\lambda},\frac{\lambda}{\lambda-1}\}$ (see proof of 
Proposition~\ref{fam}).
Explicitly, 
\[\varphi_1(X:Y:Z)=(Z:\frac{1}{\sqrt[4]{\lambda}}\;Y:X)\]
induces an isomorphism $\varphi_1:W_\lambda\to W_{\frac{1}{\lambda}}$, and
\[\varphi_2(X:Y:Z)=(Z-X:\zeta_8Y:Z)\]
induces an isomorphism $\varphi_2:W_\lambda\to W_{1-\lambda}$; here $\zeta_8$
is a primitive eighth root of 1. The other
isomorphisms can be obtained from these two by composition, since
$s_1(\lambda)=\frac{1}{\lambda}$ and $s_2(\lambda)=1-\lambda$ generate the
above group of six transformations.\\
The following relations can be checked by straightforward calculation:
\begin{prop}
\begin{eqnarray*}
\varphi_1^{-1}\circ\sigma_{\frac{1}{\lambda}}\circ\varphi_1 = (-\rho)_\lambda
&& \varphi_2^{-1}\circ\sigma_{1-\lambda}\circ\varphi_2 = (-\tau)_\lambda\\
\varphi_1^{-1}\circ\rho_{\frac{1}{\lambda}}\circ\varphi_1 = (-\sigma)_\lambda
&& \varphi_2^{-1}\circ\rho_{1-\lambda}\circ\varphi_2 = \rho_\lambda\\
\varphi_1^{-1}\circ\tau_{\frac{1}{\lambda}}\circ\varphi_1 = \tau_\lambda\phantom{(-),}
&& \varphi_2^{-1}\circ\tau_{1-\lambda}\circ\varphi_2 = (-\sigma)_\lambda.\\
\end{eqnarray*}
\end{prop}
In these equations the sign of the automorphism on the right hand side depends
on the choice of the square roots $\zeta(\lambda)$, $\xi(\lambda)$ and
$\eta(\lambda)$. E.\,g.\ the first equation is correct in this form if we choose $\xi(\lambda)=\zeta(\frac{1}{\lambda})$.
\begin{cor}
There is an isomorphism $\bar\varphi_1:E_{-1}\to E_{-1}$ such that
$\kappa_\sigma^{(\frac{1}{\lambda})}\circ\varphi_1 =
\bar\varphi_1\circ\kappa_{-\rho}^{(\lambda)}$, i.\,e.\ $\bar\varphi_1$ makes the
following diagram commutative:
\begin{center}
$\xymatrix@=6ex{\ar@(ul,dl)_{-\rho_\lambda}\ar[r]^{\varphi_1}\ar[d]_{\kappa_{-\rho}} W_\lambda&\ar@(ur,dr)^{\sigma_\frac{1}{\lambda}}\ar[d]^{\kappa_\sigma}W_{\frac{1}{\lambda}}\\
\ar[r]^{\bar\varphi_1}E_{-1}&E_{-1}}
$
\end{center}
Explicitly, $\bar\varphi_1$ is the translation by $(0:0:1)$ .\\[2mm]
Such an isomorphism exists in an analogous way for each of the isomorphisms
$W_\lambda\to W_{\lambda'}$, $\lambda'\in\{\lambda,\frac{1}{\lambda}, 1~-~\lambda,$ 
$\frac{1}{1-\lambda}, 1~-~\frac{1}{\lambda}, \frac{\lambda}{\lambda-1}\}$.
\end{cor}
\subsection{Ramification points}
\label{ramif}
In this section we calculate the ramification points and the critical points
of the two-sheeted coverings $W_\lambda\to E_{-1}$, i.\,e.\ the fixed points
of the involutions $\pm\sigma$, $\pm\rho$ and $\pm\tau$ on $W_\lambda$ and
their images on $E_{-1}$ under $\kappa_\sigma$ (resp.\ $\kappa_{-\sigma}$,
$\kappa_{\pm\rho}$ and $\kappa_{\pm\tau}$).

\begin{prop}
For any $\alpha\in\{\pm\sigma,\pm\rho,\pm\tau\}$ the fixed points of
$\alpha$ on $W_\lambda$ are
\[P^\nu_\alpha:=(X_\alpha:i^\nu:Z_\alpha),\ \nu=0,1,2,3,\]
where $X_\alpha$ and $Z_\alpha$ are defined as follows using the notation $\zeta$, 
$\xi$ and $\eta$ as in the previous section:
\begin{eqnarray*}
\begin{array}{lcllcl}
X_\sigma&=&(1+\zeta)Z_\sigma, \quad\quad&
Z_\sigma&=&\frac{1}{\sqrt{\zeta(1+\zeta)}},\\
X_\rho&=&\sqrt[4]{\frac{1+\xi}{\xi^2(1-\xi)}},&
Z_\rho&=&(1-\xi)X_\rho\ \qquad  \mbox{and}\\
X_\tau&=&\eta Z_\tau,&
Z_\tau&=&\frac{1}{\sqrt{i\eta(\eta-1)}}.
\end{array}
\end{eqnarray*}
As in section \ref{others} the fixed points of $-\sigma$ (resp.\ $-\rho$, $-\tau$)
are obtained from these formulas by replacing $\zeta$ (resp.\ $\xi$, $\eta$)
by its negative.\\ A different choice of a root in these expressions results
in a permutation of the four points $P_\alpha^\nu$, $\nu=0,1,2,3$.
\end{prop}
\begin{proof}
Recall that
$W_\lambda\cap\{Y=0\}=\{P_0,P_1,P_\infty,P_\lambda\}$, and none of these
points is fixed by any of the six involutions. So we may assume all fixed
points in the form $(X:1:Z)$. If such a point is fixed by $\sigma$, we see
from Proposition~\ref{formeln} that $Z=\frac{1}{\zeta}(X-Z)$ or
equivalently $X=Z(1+\zeta)$. Inserting this into the equation for $W_\lambda$
we find
\[1=Z^4(1+\zeta)\zeta(1+\zeta-\lambda)=Z^4\zeta^2(1+\zeta)^2,\]
which gives $Z_\sigma$. The calculations for $\rho$ and
$\tau$ are similar.
\end{proof}
Remark that $P^0_\alpha$, $P^1_\alpha$, $P^2_\alpha$, $P^3_\alpha$ form an orbit under the automorphism
$\drehung$. This can also be seen from the description of the automorphisms in Section
\ref{G}.\\[2mm]
The critical points of $\kappa_\alpha$ are the images of the $P_\alpha^\nu$,
e.\,g.: 
\begin{eqnarray*}
\kappa_\sigma(P_\sigma^\nu)&=&\left((-1)^{\nu+1}(1-\zeta):i^\nu\sqrt{1-\zeta}\frac{2\zeta}{\sqrt{\zeta(1+\zeta)}}:(1+\zeta)\frac{1+\zeta-\lambda}{\zeta(1+\zeta)}\right)\\
&=&\left((-1)^{\nu+1}(1-\zeta):i^\nu\sqrt{\frac{1-\zeta}{1+\zeta}}2\sqrt{\zeta}:1+\zeta\right)\\
&=&\left((-1)^{\nu+1}\lambda:2i^\nu\sqrt{\zeta\lambda}:(1+\zeta)^2\right)
\end{eqnarray*}
Proceeding in the same way for the other five involutions we obtain
\begin{cor} \label{rampointscor}
The images in $E_{-1}$ of the fixed points of $\pm\sigma$, $\pm\rho$ and
$\pm\tau$ under $\kappa_{\pm\sigma}$, $\kappa_{\pm\rho}$ and $\kappa_{\pm\tau}$
resp., are
\begin{eqnarray*}
Q_\sigma&=&\left(-\lambda:2\sqrt{\zeta\lambda}:(1+\zeta)^2\right)\\
Q_{-\sigma}&=&\left(\lambda:2\sqrt{\zeta\lambda}:(1-\zeta)^2\right)\\
Q_\rho&=&\left(\xi-1:2\eta(1-\xi)\sqrt{\xi}:1+\xi\right)\\
Q_{-\rho}&=&\left(\xi+1:2\eta(1+\xi)\sqrt{\xi}:1-\xi\right)\\
Q_\tau&=&\left(1-\eta:2\sqrt{\eta}\sqrt{\frac{\eta-1}{\eta+1}}:\eta+1\right)\\
Q_{-\tau}&=&\left(1+\eta:2\sqrt{\eta}\sqrt{\frac{\eta+1}{\eta-1}}:\eta-1\right)\\
\end{eqnarray*}
and their orbits under the automorphism $\bar\drehung$.
\end{cor}


\section{Intersecting origami curves} \label{sectiondrei}
Now, finally we show (in \ref{thm}) that in the moduli space $M_3$ the 
Teichm\"uller curve $C(W)$ of the quaternion origami intersects infinitely 
many other origami curves. The combinatorial description of these origamis 
is given in \ref{dq}.

\subsection{The main result}
\label{thm}
Recall from Proposition~\ref{kprop} that we have for each $\lambda \in \PPP$ the degree 2 
morphism $\kappa = \kappa_{\sigma_\lambda}: W_{\lambda} \rightarrow E_{-1}$. The 
critical points of $\kappa$ form an orbit under the order 4 automorphism 
$\bar{\drehung}$ and are given explicitly in Corollary \ref{rampointscor}. These are 
the main ingredients for the following theorem.
\begin{thm}
\label{thethm}
The origami curve $C(W)$ intersects infinitely many other origami curves.
\end{thm}
\begin{proof}
The idea of the proof is to compose $\kappa_\sigma$ (or one of the
other morphisms $W_\lambda \to E_{-1}$) with a morphism $\varphi:E_{-1}\to
E'$ to an elliptic curve $E'$ that maps
the four ramification points $Q^\nu_\sigma=\bar\drehung^\nu(Q_\sigma),\nu=0,1,2,3$ to the same point. Then
$\varphi\circ\kappa_\sigma$ is an origami.\\[1mm]
{\sc Claim.} Such a morphism exists if and only if $Q_\sigma:=Q^0_\sigma$ is a
torsion point on $E_{-1}$ (with the usual group structure, i.\,e. with
$0:=(0:1:0)$ as origin).\\[2mm]
{\it Proof of the claim.} By composing with an isogeny $E'\to E_{-1}$ if
necessary we may assume that $\varphi$ is the product of a
translation with multiplication by $n$ on $E_{-1}$ (for some $n\ge2$). Hence
$\varphi(Q^\nu_\sigma)=\varphi(Q_\sigma)$ for $\nu=0,\dots,3$ if and only if
there is a point $P\in E_{-1}(\CC)$ such that $n\cdot (Q^\nu_\sigma-P)=Q'$ with
the same point $Q'$ for $\nu=0,\dots,3$. On the other hand, since
$Q^\nu_\sigma=\bar\drehung^\nu(Q_\sigma)$, we have
$n\cdot\bar\drehung^\nu(Q_\sigma)=n\cdot P + Q'$ for $\nu=0,\dots,3$. Since
$\bar\drehung$ commutes with multiplication by $n$, this implies in particular
that $n\cdot Q_\sigma$ is a fixed point of $\bar\drehung$, i.\,e. either $0$
or $(0:0:1)$. In both cases $2n\cdot Q_\sigma=0$. This proves the claim.\\[1mm]
From Section \ref{ramif} we know that 
\begin{eqnarray*}
Q_\sigma=Q_\sigma(\lambda)&=&\left(-\lambda:2\sqrt{\zeta\lambda}:(1+\zeta)^2\right)=\left(-\frac{\lambda}{(1+\zeta)^2}:2\frac{\sqrt{\zeta\lambda}}{(1+\zeta)^2}:1\right)\\
&=&\left(\frac{\zeta-1}{\zeta+1}:2\frac{\sqrt{\zeta\lambda}}{(1+\zeta)^2}:1\right)  
\in E_{-1}: y^2 = x^3 - x
\end{eqnarray*}
As $W_\lambda$ runs through $C(W)$, $\lambda$ takes on every value in
$\PPP = \PP^1(\CC)-\{0,1,\infty\}$. Since $\zeta^2=1-\lambda$, $\zeta$ omits only the
values $0$, $\infty$ and $\pm 1$, and the same holds for
$\frac{\zeta-1}{\zeta+1}$. This shows that for every $Q\in E_{-1}$ which is
not a 2-torsion point we find $\lambda\in\PPP$ such that
$Q_\sigma(\lambda)=Q$. Together with the claim this proves that for every
torsion point $T$ of $E_{-1}$ of (exact) order 
$n\ge3$ there is $\lambda_T\in\PPP$ such that
$[n]\circ\kappa_\sigma:W_{\lambda_T}\to E_{-1}$ defines an origami $D_T$ which is
different from $W$. By construction $C(D_T)$ and $C(W)$ intersect in
$W_{\lambda_T}$.
\end{proof}

\subsection{Nice double coverings of the trivial origami} \label{dq}

In this section we describe the origamis intersecting $C(W)$.\\[2mm]
Recall from the last section that they are given as composition of the 
double cover $\kappa: W_\lambda \rightarrow 
E_{-1}$ and $[n]: E_{-1} \rightarrow E_{-1}$ , the multiplication
with $n$. The four critical points of $\kappa$ are $n$-torsion 
points but not $2$-torsion points. They form an orbit under the 
automorphism $\bar{\drehung}$ on $E_{-1}$ which has order four and 
fixed points $0$ and $(0:0:1)$. Furthermore, 
$\bar{\drehung}$ is induced by the automorphism $c$ on $W_{\lambda}$
and $c$ as well as $c^2$ have four fixed points.\\[2mm]   
Call $\tr$ the origami with the translation structure 
defined by the unit lattice on $E_{-1}$. Furthermore, call
$\trn$, respectively $D$, the origami given as the pullback via
$[n]$ on $E_{-1}$, respectively via $[n] \circ \kappa$ on $W_{\lambda}$.
Then $D$ is the desired origami and $\trn$ consists of a large 
square of length $n$ subdivided into $n^2$ squares 
(see Figure \ref{0qn}) where opposite sides are glued.
We will show that the properties listed in the last
paragraph uniquely determine $D$.\\[2mm]
Remark that the automorphisms $\bar{\drehung}$ and $\bar{\drehung}^{-1}$
are affine with derivative 
\[ S := \begin{pmatrix}0&-1\\1&0\end{pmatrix} \mbox{\; and \;} 
   S^{-1} := \begin{pmatrix}0&1\\-1&0\end{pmatrix}, \mbox{ resp.,} \]
since they are the only automorphisms of order four fixing $0$.
The lift $c$ is thus also affine with the same derivative.

\begin{prop} \label{schnittorigamis}
Let $P$ be a $n$-torsion point of $\trn$ which is not a $2$-torsion point.
Then, there is a uniquely determined origami $D = (p: X \rightarrow \equ)$ that has a 
translation covering $q: X \rightarrow \trn$ of degree 2 with $p = [n]\circ q$ such that 
the following properties hold:
\begin{enumerate}
\item \label{cond1} There is an affine automorphism $\drehung$ on $X$ with derivative 
$D(\drehung) = S$
that descends to the rotation $\abar$ on $\trn$ with derivative $S$ and center $0$. 
The automorphism $\drehung$ and its square $\drehung^2$ have four fixed points on $X$.
\item \label{cond2} The critical points of $q$ are $P$, 
$Q := \abar(P)$, $R := \abar^2(P)$ and $S := \abar^3(P)$.
\end{enumerate}
\end{prop}

\begin{center}
\setlength{\unitlength}{.5cm}
\begin{minipage}{4cm}

\begin{picture}(7,7)
\put(1,1){\framebox(1,1){}}
\put(1,2){\framebox(1,1){}}
\put(1,3){\framebox(1,1){}}
\put(1,4){\framebox(1,1){}}
\put(1,5){\framebox(1,1){}}
\put(2,1){\framebox(1,1){}}
\put(2,2){\framebox(1,1){}}
\put(2,3){\framebox(1,1){}}
\put(2,4){\framebox(1,1){}}
\put(2,5){\framebox(1,1){}}
\put(2,1){\framebox(1,1){}}
\put(2,2){\framebox(1,1){}}
\put(2,3){\framebox(1,1){}}
\put(2,4){\framebox(1,1){}}
\put(2,5){\framebox(1,1){}}
\put(3,1){\framebox(1,1){}}
\put(3,2){\framebox(1,1){}}
\put(3,3){\framebox(1,1){}}
\put(3,4){\framebox(1,1){}}
\put(3,5){\framebox(1,1){}}
\put(4,1){\framebox(1,1){}}
\put(4,2){\framebox(1,1){}}
\put(4,3){\framebox(1,1){}}
\put(4,4){\framebox(1,1){}}
\put(4,5){\framebox(1,1){}}
\put(5,1){\framebox(1,1){}}
\put(5,2){\framebox(1,1){}}
\put(5,3){\framebox(1,1){}}
\put(5,4){\framebox(1,1){}}
\put(5,5){\framebox(1,1){}}

\multiput(0,2.03)(.5,0){14}{\line(1,0){.3}}
\multiput(0,5.03)(.5,0){14}{\line(1,0){.3}}
\put(6.4,2.1){$h_1$}
\put(6.4,5.1){$h_2$}
\multiput(1.97,0)(0,.5){14}{\line(0,1){.3}}
\multiput(4.97,0)(0,.5){14}{\line(0,1){.3}}
\put(5.2,6.4){$v_1$}
\put(2.2,6.4){$v_2$}

\put(1.4,1.4){$A$}
\put(5.1,1.4){$B$}
\put(5.1,5.1){$C$}
\put(1.4,5.1){$D$}

\put(2,2){\circle*{.1}}
\put(2,2){\circle{.2}}
\put(5,2){\circle*{.1}}
\put(5,2){\circle{.2}}
\put(5,5){\circle*{.1}}
\put(5,5){\circle{.2}}
\put(2,5){\circle*{.1}}
\put(2,5){\circle{.2}}

\put(3.2,1.4){$P$}
\put(5.2,2.4){$Q$}
\put(3.4,5.2){$R$}
\put(1.4,3.4){$S$}

\put(3,2){\circle*{.25}}
\put(5,3){\circle*{.25}}
\put(4,5){\circle*{.25}}
\put(2,4){\circle*{.25}}

\put(0.65,1){$0$}
\put(1,1){\circle*{.1}}
\put(1,1){\circle{.2}}

\put(3.2,3.7){$M$}
\put(3.5,3.5){\circle*{.1}}
\put(3.5,3.5){\circle{.2}}

\put(3.4,0.4){$F$}
\put(3.5,1){\circle*{.1}}
\put(3.5,1){\circle{.2}}

\put(0.1,3.4){$F'$}
\put(1,3.5){\circle*{.1}}
\put(1,3.5){\circle{.2}}

\end{picture}
\begin{center}
\captionfigure \label{0qn}
\end{center}
\end{minipage}
\hspace{2cm}
\setlength{\unitlength}{.5cm}
\begin{minipage}{4cm}
\begin{picture}(7,7)
\put(1,1){\framebox(5,5){}}

\put(3,1){\line(0,1){.3}}
\put(4,1){\line(0,1){.3}}
\put(5,1){\line(0,1){.3}}
\put(1,3.5){\line(1,0){5}}
\put(3.5,1){\line(0,1){5}}
\put(1,1){\line(1,1){5}}
\put(1,6){\line(1,-1){5}}





\put(0.65,1){$0$}
\put(1,1){\circle*{.1}}
\put(1,1){\circle{.2}}

\put(6.25,1){$0$}
\put(6,1){\circle*{.1}}
\put(6,1){\circle{.2}}

\put(6.25,6){$0$}
\put(6,6){\circle*{.1}}
\put(6,6){\circle{.2}}

\put(0.65,6){$0$}
\put(1,6){\circle*{.1}}
\put(1,6){\circle{.2}}

\put(3.6,3.6){$M$}
\put(3.5,3.5){\circle*{.1}}
\put(3.5,3.5){\circle{.2}}

\put(3.4,0.4){$F$}
\put(3.5,1){\circle*{.1}}
\put(3.5,1){\circle{.2}}

\put(3.4,6.2){$F$}
\put(3.5,6){\circle*{.1}}
\put(3.5,6){\circle{.2}}

\put(6.25,3.3){$F'$}
\put(6,3.5){\circle*{.1}}
\put(6,3.5){\circle{.2}}

\put(0.4,3.3){$F'$}
\put(1,3.5){\circle*{.1}}
\put(1,3.5){\circle{.2}}

\put(2.5,0.4){\vector(0,1){1.2}}
\put(4.5,0.4){\vector(0,1){1.2}}
\put(6.6,2.5){\vector(-1,0){1}}
\put(6.6,4.5){\vector(-1,0){1}}

\put(2.5,6.6){\vector(0,-1){1.2}}
\put(4.5,6.6){\vector(0,-1){1.2}}
\put(.6,2.4){\vector(1,0){1.2}}
\put(.6,4.4){\vector(1,0){1.2}}

\put(1.5,0.1){$\Delta_1$}
\put(4.5,0.1){$\Delta_2$}
\put(6.6,2.5){$\Delta_3$}
\put(6.6,4.5){$\Delta_4$}

\put(1.5,6.2){$\Delta_6$}
\put(4.6,6.2){$\Delta_5$}
\put(0,2.6){$\Delta_8$}
\put(0,4.6){$\Delta_7$}
\end{picture}
\begin{center}
\captionfigure \label{1qn}
\end{center}
\end{minipage}
\end{center}

\begin{proof}
First we show that there is at most one origami like this. 
Suppose $D = (p: X \rightarrow \equ)$ is an origami with a covering 
$q:X \rightarrow \trn$ fulfilling the properties in the proposition.
Draw - as shown in Figure \ref{0qn} - the horizontal lines $h_1$ 
through $P$ and $h_2$ through $R$ and the vertical lines $v_1$ 
through $Q$ and $v_2$ through $S$. Define $A:= h_1 \cap v_2$, 
$B:= h_1 \cap v_1$, $C:= h_2 \cap v_1$, $D:= h_2 \cap v_2$. Let 
$M\neq 0$ be the other fixed point of $\abar$ and let  
$F$ and $F'$ be the other $2$-torsion points of $\trn$ on the 
horizontal and the vertical geodesic starting in $0$.\\
$\trn$ is divided by the horizontal and vertical lines through $0$ and 
$M$ and the two diagonals into eight (euclidian) geodesic triangles 
$\Delta_1, \ldots, \Delta_8$ (see Figure \ref{1qn}).
We distinguish five cases: 
\begin{enumerate}
\item \label{case1} $P$ lies in the open triangle $\Delta_1, \Delta_3, \Delta_5$ or $\Delta_7$.
\item \label{case2} $P$ lies in the open triangle $\Delta_2, \Delta_4, \Delta_6$ or $\Delta_8$.
\item \label{case3} $P$ lies on one of the two diagonals.
\item \label{case4} $P$ lies on the horizontal or the vertical line 
through $0$.
\item \label{case5} $P$ lies on the horizontal or the vertical line 
through $M$.
\end{enumerate}
{\bf Case \ref{case1}:}\\


\begin{center}
\setlength{\unitlength}{.7cm}
\begin{picture}(7,7)
\put(1,1){\framebox(1,1){}}
\put(1,2){\framebox(1,1){}}
\put(1,3){\framebox(1,1){}}
\put(1,4){\framebox(1,1){}}
\put(1,5){\framebox(1,1){}}
\put(2,1){\framebox(1,1){}}
\put(2,2){\framebox(1,1){}}
\put(2,3){\framebox(1,1){}}
\put(2,4){\framebox(1,1){}}
\put(2,5){\framebox(1,1){}}
\put(2,1){\framebox(1,1){}}
\put(2,2){\framebox(1,1){}}
\put(2,3){\framebox(1,1){}}
\put(2,4){\framebox(1,1){}}
\put(2,5){\framebox(1,1){}}
\put(3,1){\framebox(1,1){}}
\put(3,2){\framebox(1,1){}}
\put(3,3){\framebox(1,1){}}
\put(3,4){\framebox(1,1){}}
\put(3,5){\framebox(1,1){}}
\put(4,1){\framebox(1,1){}}
\put(4,2){\framebox(1,1){}}
\put(4,3){\framebox(1,1){}}
\put(4,4){\framebox(1,1){}}
\put(4,5){\framebox(1,1){}}
\put(5,1){\framebox(1,1){}}
\put(5,2){\framebox(1,1){}}
\put(5,3){\framebox(1,1){}}
\put(5,4){\framebox(1,1){}}
\put(5,5){\framebox(1,1){}}

\multiput(0,2.03)(.5,0){14}{\line(1,0){.3}}
\multiput(0,5.03)(.5,0){14}{\line(1,0){.3}}
\put(6.2,2.1){$h_1$}
\put(6.2,5.1){$h_2$}
\multiput(1.97,0)(0,.5){14}{\line(0,1){.3}}
\multiput(4.97,0)(0,.5){14}{\line(0,1){.3}}
\put(5.2,6.1){$v_1$}
\put(2.2,6.1){$v_2$}

\put(1.65,1.6){$A$}
\put(5.1,1.6){$B$}
\put(5.1,5.1){$C$}
\put(1.5,5.1){$D$}

\put(2,2){\circle*{.1}}
\put(2,2){\circle{.2}}
\put(5,2){\circle*{.1}}
\put(5,2){\circle{.2}}
\put(5,5){\circle*{.1}}
\put(5,5){\circle{.2}}
\put(2,5){\circle*{.1}}
\put(2,5){\circle{.2}}

\put(3.2,1.6){$P$}
\put(5.2,2.65){$Q$}
\put(3.6,5.2){$R$}
\put(1.65,3.65){$S$}

\put(3,2){\circle*{.25}}
\put(5,3){\circle*{.25}}
\put(4,5){\circle*{.25}}
\put(2,4){\circle*{.25}}

\put(1,2){\rule[-0.06\unitlength]{0.9\unitlength}{0.12\unitlength}}
\put(2.1,2){\rule[-0.06\unitlength]{0.8\unitlength}{0.12\unitlength}}
\put(3.1,2){\rule[-0.06\unitlength]{1.8\unitlength}{0.12\unitlength}}
\put(5.1,2){\rule[-0.06\unitlength]{.9\unitlength}{0.12\unitlength}}

\put(4.96,3){\rule{0.12\unitlength}{1.9\unitlength}}
\put(4.96,5.1){\rule{0.12\unitlength}{.9\unitlength}}
\put(4.96,1){\rule{0.12\unitlength}{.9\unitlength}}
\put(4.1,5){\rule[-0.06\unitlength]{0.8\unitlength}{0.12\unitlength}}
\put(1.96,2.1){\rule{0.12\unitlength}{1.8\unitlength}}
\put(1.96,4.1){\rule{0.12\unitlength}{1.9\unitlength}}
\put(1.96,2.1){\rule{0.12\unitlength}{1.9\unitlength}}
\put(2,1){\rule{0.12\unitlength}{.9\unitlength}}

\put(2.4,1.65){$e_1$}
\put(4.2,1.65){$e_2$}
\put(5.6,1.65){$e_3$}
\put(5.1,3.7){$e_5$}
\put(5.1,5.6){$e_6$}
\put(5.1,1.2){$e_6$}
\put(4.3,4.65){$e_7$}
\put(1.3,4.4){$e_{10}$}
\put(1.3,5.7){$e_{10}$}
\put(1.35,1.25){$e_{10}$}
\put(1.3,3.2){$e_{11}$}

\put(0.65,1){$0$}
\put(1,1){\circle*{.1}}
\put(1,1){\circle{.2}}

\put(3.6,3.6){$M$}
\put(3.5,3.5){\circle*{.1}}
\put(3.5,3.5){\circle{.2}}

\put(3.4,0.4){$F$}
\put(3.5,1){\circle*{.1}}
\put(3.5,1){\circle{.2}}

\put(0.4,3.4){$F$}
\put(1,3.5){\circle*{.1}}
\put(1,3.5){\circle{.2}}

\end{picture}
\begin{center}
\captionfigure \label{2case1}
\end{center}
\end{center}
Suppose $P$ lies inside $\Delta_1$. If it lies in one of the other three triangles 
one can interchange the roles of $P$, $Q$, $R$ and $S$.\\
Since $P$ lies inside the open triangle $\Delta_1$,  the lines $h_1,h_2,v_1,v_2$ are 
pairwise different, the sets $\{A,B,C,D\}$ and $\{P,Q,R,S\}$ are disjoint and  we have: 
$P$ lies between $A$ and $B$, $Q$ lies between $B$ and $C$, $R$ lies between $D$ and $C$, $S$ lies between $A$ and $D$ \footnote{\label{orientation} Orientation is always left to right for horizontal directions and bottom to top for vertical directions.} 
(compare Figure \ref{2case1}). Define the segments $e_1 := AP, e_2 := PB$, $e_3 := BA$, $e_5 := QC$, $e_6 := CB$, $e_7 := RC$, $e_{10} := SA$, $e_{11} := AS$ (compare Footnote \ref{orientation}). Removing $e_1, e_2, e_3, e_5, e_6, e_7, e_{10}, e_{11}$ from $\trn$ gives a simply connected open region $U$. 
Since the critical points $P$, $Q$, $R$ and $S$ are contained in the border of 
$U$, one can get any degree 2 cover of $\trn$ which has those as critical points by 
taking two copies $U_1$ and $U_2$ of $U$ and glueing the borders of $\bar{U}_1$ and $\bar{U}_2$.\\
Now, take a glueing that gives $q:X \rightarrow \trn$. We will show that it is uniquely determined by the properties required in the proposition.\\
Let's label all horizontal and vertical edges $e$ of length 1 between two neighbouring $n$-torsion points of $\trn$ by a label $l(e) \in \ZZ/2\ZZ$ in the following way:
\begin{itemize}
\item[-] $l(e) = \bar{1}$, if there is a change of the two leaves at $e$, 
i.e. if for a preimage $Z$ of a point in the interior of $e$, 
each neigbourhood of $Z$ intersects $U_1$ and $U_2$.
\item[-] $l(e) = \bar{0}$, if there is no change of the leaves at $e$, 
i.e. if for a preimage $Z$ of a point in the interior of $e$, 
there is a neigbourhood of $Z$ that lies completely inside $U_1$ or $U_2$.
\end{itemize}
By construction we have $l(e) = \bar{1} \Rightarrow e \in  \{e_1, e_2, e_3, e_5, e_6, e_7, e_{10}, e_{11}\}$.\\
Furthermore, an $n$-torsion point on $\trn$ that is intersection of the four edges $e, e', e'', e'''$ is 
\begin{itemize}
\item[-] ramified $\Leftrightarrow$ $l(e) + l(e') + l(e'')+ l(e''') = \bar{1}$. \hspace{\fill} (*)
\item[-] unramified $\Leftrightarrow$ $l(e) + l(e') + l(e'')+ l(e''') = \bar{0}$.\hspace{\fill} (**)
\end{itemize}
Hence  we get by Condition \ref{cond2} in the proposition:
\begin{eqnarray}
\label{vzwgl1}
&&l(e_1) + l(e_2) = 1, \quad\quad l(e_5) = 1,\\ \nonumber
&&l(e_7) = 1, l(e_{10}) + l(e_{11}) = 1  \quad\quad\quad\quad\quad
\mbox{ ($P$, $Q$, $R$, $S$ ramified)}\\
\label{vzwgl2}
&&l(e_1) + l(e_{11}) + l(e_3) + l(e_{10}) = 0, \quad\quad
l(e_2) + l(e_{6}) + l(e_3) = 0,\\ \nonumber
&&l(e_5) + l(e_{6}) + l(e_7)  = 0  \quad\quad\quad\quad\quad\;\;
\mbox{ ($A$, $B$, $C$, $D$ unramified)}
\end{eqnarray}
By Condition \ref{cond1} in the proposition the automorphism 
$\drehung$ has four fixed points on $X$. Since $\abar$ on $\trn$ 
has the two fixed points $0$ and $M$, the four fixed points of $\drehung$ 
are their preimages $0_1 := q^{-1}(\{0\}) \cap U_1$, $0_2 := q^{-1}(\{0\}) 
\cap U_2$, $M_1 := q^{-1}(\{M\}) \cap U_1$ and $M_2 := q^{-1}(\{M\})
\cap U_2$.\\
Obviously, they are also fixed points of $\drehung^2$ and by condition \ref{cond1} 
they are all of them. From this it follows that $l(e_2) = 1$: If not the vertical 
geodesic starting in $M_1$ would be a closed geodesic of length $n$ and 
the $2$-torsion point $F$ of $\trn$ on it would be a fixed point of $\drehung^2$.\\ 
Similarly, we get $l(e_5) + l(e_{11}) = 1$.
Hence we have
\begin{eqnarray}
\label{autogl}
l(e_2) = 1 \quad \mbox{ and } \quad l(e_5) + l(e_{11}) = 1.
\end{eqnarray}
From the equations (\ref{vzwgl1}), (\ref{vzwgl2}) and (\ref{autogl}) it follows that:
\begin{eqnarray*}
\begin{array}{c@{\;\;,\quad}c@{\;\;,\quad}c@{\;\;,\quad}c}
l(e_1) = 0 & l(e_2) = 1 & l(e_3) = 1 & l(e_5) = 1,\\
l(e_6) =  0& l(e_7) = 1 & l(e_{10}) = 1 & l(e_{11})  = 0
\end{array}
\end{eqnarray*}
Hence, the origami $D$ is uniquely determined by the conditions in the proposition: 
the leaves are switched at the edges $e_2$, $e_3$, $e_5$, $e_7$ and $e_{10}$. \\[3mm]
Conversely, this origami fulfills the condition of the proposition:\\
In order to see that the automorphism $\abar$ of $\trn$ can be lifted 
to $X$ we consider the fundamental group $H$ of 
$\trn \backslash \{P,Q,R,S\}$. It is a free group in the generators 
$\alpha$, $\beta$, $\gamma_1$, $\gamma_2$, $\gamma_3$, where 
$\alpha$ and $\beta$ are the closed horizontal 
and vertical geodesic lines through $M$ and $\gamma_1$, $\gamma_2$, 
$\gamma_3$ are the 
positively oriented simple loops around $P$, $Q$ and $R$. The loop 
$\gamma_4$ around $S$ is equivalent to $\gamma_3^{-1}\gamma_2^{-1}\gamma_1^{-1}$.\\
By the definition of the glueings the monodromy $H \rightarrow \ZZ/2\ZZ$ is given by 
$g \mapsto 1$ for each generator $g$. The fundamental group $U$ of 
$X \backslash q^{-1}(\{P,Q,R,S\})$ is the kernel of the monodromy. 
Thus it is the index $2$ subgroup of words of even length.\\
The automorphism of $H$ induced by $\abar$  is given by
\[\alpha \mapsto \beta,\quad \beta \mapsto \alpha^{-1},\quad
  \gamma_1\mapsto \gamma_2,\quad \gamma_2 \mapsto \gamma_3,\quad
  \gamma_3 \mapsto  \gamma_4 = \gamma_3^{-1}\gamma_2^{-1}\gamma_1^{-1}.\]
It preserves the parity of the length of words and thus it preserves $U$. 
It follows that
  $\abar$ can be lifted to $\drehung$ on $X$.
One checks that $\drehung$ and $\drehung^2$ have the fixed point set 
$q^{-1}(\{0,M\})$. Thus Condition \ref{cond1} is fulfilled.\\
Condition \ref{cond2} is fulfilled since (*) holds for $P$, $Q$, $R$ and 
$S$ and (**) holds for $A$,$B$,$C$ and $D$. This finishes the 
proof for Case \ref{case1}.\\[3mm]
The other four cases work in the same way.
\begin{center}
\setlength{\unitlength}{.6cm}
\begin{minipage}{4cm}
{\bf Case \ref{case2}}:\\
\begin{picture}(7,7)
\label{bildcase2}
\put(1,1){\framebox(1,1){}}
\put(1,2){\framebox(1,1){}}
\put(1,3){\framebox(1,1){}}
\put(1,4){\framebox(1,1){}}
\put(1,5){\framebox(1,1){}}
\put(2,1){\framebox(1,1){}}
\put(2,2){\framebox(1,1){}}
\put(2,3){\framebox(1,1){}}
\put(2,4){\framebox(1,1){}}
\put(2,5){\framebox(1,1){}}
\put(2,1){\framebox(1,1){}}
\put(2,2){\framebox(1,1){}}
\put(2,3){\framebox(1,1){}}
\put(2,4){\framebox(1,1){}}
\put(2,5){\framebox(1,1){}}
\put(3,1){\framebox(1,1){}}
\put(3,2){\framebox(1,1){}}
\put(3,3){\framebox(1,1){}}
\put(3,4){\framebox(1,1){}}
\put(3,5){\framebox(1,1){}}
\put(4,1){\framebox(1,1){}}
\put(4,2){\framebox(1,1){}}
\put(4,3){\framebox(1,1){}}
\put(4,4){\framebox(1,1){}}
\put(4,5){\framebox(1,1){}}
\put(5,1){\framebox(1,1){}}
\put(5,2){\framebox(1,1){}}
\put(5,3){\framebox(1,1){}}
\put(5,4){\framebox(1,1){}}
\put(5,5){\framebox(1,1){}}

\multiput(0,2.03)(.5,0){14}{\line(1,0){.3}}
\multiput(0,5.03)(.5,0){14}{\line(1,0){.3}}
\put(6.3,2.1){$h_1$}
\put(6.3,5.1){$h_2$}
\multiput(1.97,0)(0,.5){14}{\line(0,1){.3}}
\multiput(4.97,0)(0,.5){14}{\line(0,1){.3}}
\put(5.2,6.3){$v_1$}
\put(2.2,6.3){$v_2$}

\put(2.1,2.1){$A$}
\put(4.6,2.1){$B$}
\put(5.1,5.1){$C$}
\put(1.5,5.05){$D$}

\put(2,2){\circle*{.1}}
\put(2,2){\circle{.2}}
\put(5,2){\circle*{.1}}
\put(5,2){\circle{.2}}
\put(5,5){\circle*{.1}}
\put(5,5){\circle{.2}}
\put(2,5){\circle*{.1}}
\put(2,5){\circle{.2}}

\put(3.55,1.5){$P$}
\put(4,2){\circle*{.25}}
\put(4.4,4.2){$Q$}
\put(5,4){\circle*{.25}}
\put(3.15,5.15){$R$}
\put(3,5){\circle*{.25}}
\put(1.65,3.25){$S$}
\put(2,3){\circle*{.25}}

\put(1,2){\rule[-0.06\unitlength]{0.9\unitlength}{0.12\unitlength}}
\put(2.1,2){\rule[-0.06\unitlength]{1.8\unitlength}{0.12\unitlength}}
\put(4.1,2){\rule[-0.06\unitlength]{.8\unitlength}{0.12\unitlength}}
\put(5.1,2){\rule[-0.06\unitlength]{.9\unitlength}{0.12\unitlength}}
\put(4.96,4.1){\rule{0.12\unitlength}{.8\unitlength}}
\put(4.96,5.1){\rule{0.12\unitlength}{.9\unitlength}}
\put(4.96,1){\rule{0.12\unitlength}{.9\unitlength}}
\put(3.1,5){\rule[-0.06\unitlength]{1.8\unitlength}{0.12\unitlength}}
\put(1.96,3.1){\rule{0.12\unitlength}{2.9\unitlength}}
\put(1.96,1){\rule{0.12\unitlength}{.9\unitlength}}
\put(1.96,2.15){\rule{0.12\unitlength}{.75\unitlength}}


\put(2.4,1.65){$e_1$}
\put(4.2,1.65){$e_2$}
\put(5.55,1.65){$e_3$}
\put(1,1.65){$e_3$}
\put(5.1,4.3){$e_5$}
\put(5.1,5.65){$e_6$}
\put(5.1,1.15){$e_6$}
\put(4.3,5.1){$e_7$}
\put(1.25,4.4){$e_{10}$}
\put(1.25,5.65){$e_{10}$}
\put(1.25,1.15){$e_{10}$}
\put(1.25,2.4){$e_{11}$}

\put(0.65,1){$0$}
\put(1,1){\circle*{.1}}
\put(1,1){\circle{.2}}

\put(3.6,3.6){$M$}
\put(3.5,3.5){\circle*{.1}}
\put(3.5,3.5){\circle{.2}}

\put(3.4,0.4){$F$}
\put(3.5,1){\circle*{.1}}
\put(3.5,1){\circle{.2}}

\put(0.2,3.3){$F'$}
\put(1,3.5){\circle*{.1}}
\put(1,3.5){\circle{.2}}

\end{picture}
\begin{center}
\captionfigure \label{3case2}
\end{center}
\end{minipage}
\hspace{2cm}
\setlength{\unitlength}{.6cm}
\begin{minipage}{4cm}
{\bf Case \ref{case3}}:\\
\begin{picture}(7,7)
\label{bildcase3}
\put(1,1){\framebox(1,1){}}
\put(1,2){\framebox(1,1){}}
\put(1,3){\framebox(1,1){}}
\put(1,4){\framebox(1,1){}}
\put(1,5){\framebox(1,1){}}
\put(2,1){\framebox(1,1){}}
\put(2,2){\framebox(1,1){}}
\put(2,3){\framebox(1,1){}}
\put(2,4){\framebox(1,1){}}
\put(2,5){\framebox(1,1){}}
\put(2,1){\framebox(1,1){}}
\put(2,2){\framebox(1,1){}}
\put(2,3){\framebox(1,1){}}
\put(2,4){\framebox(1,1){}}
\put(2,5){\framebox(1,1){}}
\put(3,1){\framebox(1,1){}}
\put(3,2){\framebox(1,1){}}
\put(3,3){\framebox(1,1){}}
\put(3,4){\framebox(1,1){}}
\put(3,5){\framebox(1,1){}}
\put(4,1){\framebox(1,1){}}
\put(4,2){\framebox(1,1){}}
\put(4,3){\framebox(1,1){}}
\put(4,4){\framebox(1,1){}}
\put(4,5){\framebox(1,1){}}
\put(5,1){\framebox(1,1){}}
\put(5,2){\framebox(1,1){}}
\put(5,3){\framebox(1,1){}}
\put(5,4){\framebox(1,1){}}
\put(5,5){\framebox(1,1){}}

\multiput(0,2.03)(.5,0){14}{\line(1,0){.3}}
\multiput(0,5.03)(.5,0){14}{\line(1,0){.3}}
\put(6.3,2.1){$h_1$}
\put(6.3,5.1){$h_2$}
\multiput(1.97,0)(0,.5){14}{\line(0,1){.3}}
\multiput(4.97,0)(0,.5){14}{\line(0,1){.3}}
\put(5.2,6.3){$v_1$}
\put(2.2,6.3){$v_2$}

\put(2.1,2.1){$P$}
\put(4.55,2.1){$Q$}
\put(5.1,5.1){$R$}
\put(1.5,5.1){$S$}

\put(2,2){\circle*{.3}}
\put(5,2){\circle*{.3}}
\put(5,5){\circle*{.3}}
\put(2,5){\circle*{.3}}

\put(2.1,2){\rule[-0.06\unitlength]{2.8\unitlength}{0.12\unitlength}}
\put(5.1,2){\rule[-0.06\unitlength]{.9\unitlength}{0.12\unitlength}}
\put(1,2){\rule[-0.06\unitlength]{.9\unitlength}{0.12\unitlength}}
\put(4.96,1){\rule{0.12\unitlength}{.9\unitlength}}
\put(4.96,5.1){\rule{0.12\unitlength}{.9\unitlength}}
\put(1.96,2.1){\rule{0.12\unitlength}{2.9\unitlength}}
\put(1.96,1){\rule{0.12\unitlength}{.9\unitlength}}
\put(1.96,5.1){\rule{0.12\unitlength}{.9\unitlength}}


\put(3.4,1.65){$e_1$}
\put(5.4,1.65){$e_2$}
\put(1.1,1.65){$e_2$}
\put(5.1,5.65){$e_5$}
\put(5.1,1.15){$e_5$}
\put(1.3,3.4){$e_{10}$}
\put(1.3,1.2){$e_{11}$}
\put(1.3,5.7){$e_{11}$}

\put(0.65,1){$0$}
\put(1,1){\circle*{.1}}
\put(1,1){\circle{.2}}

\put(3.6,3.6){$M$}
\put(3.5,3.5){\circle*{.1}}
\put(3.5,3.5){\circle{.2}}

\put(3.4,0.4){$F$}
\put(3.5,1){\circle*{.1}}
\put(3.5,1){\circle{.2}}

\put(0.2,3.3){$F'$}
\put(1,3.5){\circle*{.1}}
\put(1,3.5){\circle{.2}}

\end{picture}
\begin{center}
\captionfigure \label{4case3}
\end{center}
\end{minipage}\\
\end{center}

\begin{center}
\setlength{\unitlength}{.6cm}
\begin{minipage}{4cm}
{\bf Case \ref{case4}}:\\
\begin{picture}(7,7)
\label{bildcase3}
\put(1,1){\framebox(1,1){}}
\put(1,2){\framebox(1,1){}}
\put(1,3){\framebox(1,1){}}
\put(1,4){\framebox(1,1){}}
\put(1,5){\framebox(1,1){}}
\put(2,1){\framebox(1,1){}}
\put(2,2){\framebox(1,1){}}
\put(2,3){\framebox(1,1){}}
\put(2,4){\framebox(1,1){}}
\put(2,5){\framebox(1,1){}}
\put(2,1){\framebox(1,1){}}
\put(2,2){\framebox(1,1){}}
\put(2,3){\framebox(1,1){}}
\put(2,4){\framebox(1,1){}}
\put(2,5){\framebox(1,1){}}
\put(3,1){\framebox(1,1){}}
\put(3,2){\framebox(1,1){}}
\put(3,3){\framebox(1,1){}}
\put(3,4){\framebox(1,1){}}
\put(3,5){\framebox(1,1){}}
\put(4,1){\framebox(1,1){}}
\put(4,2){\framebox(1,1){}}
\put(4,3){\framebox(1,1){}}
\put(4,4){\framebox(1,1){}}
\put(4,5){\framebox(1,1){}}
\put(5,1){\framebox(1,1){}}
\put(5,2){\framebox(1,1){}}
\put(5,3){\framebox(1,1){}}
\put(5,4){\framebox(1,1){}}
\put(5,5){\framebox(1,1){}}

\multiput(0,1.03)(.5,0){14}{\line(1,0){.3}}
\multiput(0,6.03)(.5,0){14}{\line(1,0){.3}}
\put(6.7,1.1){$h_1$}
\put(6.7,6.1){$h_2$}
\multiput(0.97,0)(0,.5){14}{\line(0,1){.3}}
\multiput(5.97,0)(0,.5){14}{\line(0,1){.3}}
\put(5.5,6.4){$v_1$}
\put(1.1,6.4){$v_2$}

\put(2.1,0.5){$P$}
\put(2.1,6.15){$P$}
\put(6.1,2.1){$Q$}
\put(0.5,2.1){$Q$}
\put(4.6,6.1){$R$}
\put(4.6,0.5){$R$}
\put(0.4,4.6){$S$}
\put(6.1,4.6){$S$}

\put(2,1){\circle*{.3}}
\put(2,6){\circle*{.3}}
\put(6,2){\circle*{.3}}
\put(1,2){\circle*{.3}}
\put(5,6){\circle*{.3}}
\put(5,1){\circle*{.3}}
\put(1,5){\circle*{.3}}
\put(6,5){\circle*{.3}}

\put(2.1,1){\rule[-0.06\unitlength]{2.8\unitlength}{0.12\unitlength}}
\put(5.1,1){\rule[-0.06\unitlength]{.8\unitlength}{0.12\unitlength}}
\put(1.1,1){\rule[-0.06\unitlength]{.8\unitlength}{0.12\unitlength}}
\put(5.96,2.1){\rule{0.12\unitlength}{2.8\unitlength}}
\put(5.96,5.1){\rule{0.12\unitlength}{.8\unitlength}}
\put(5.96,1.1){\rule{0.12\unitlength}{.8\unitlength}}

\put(2.7,.65){$e_1$}
\put(4.1,.65){$e_1$}
\put(5.2,.65){$e_2$}
\put(1.2,.65){$e_2$}
\put(6.1,3.4){$e_4$}
\put(6.1,5.6){$e_5$}
\put(6.1,1.5){$e_5$}

\put(0.65,1){$0$}
\put(1,1){\circle*{.1}}
\put(1,1){\circle{.2}}

\put(3.6,3.6){$M$}
\put(3.5,3.5){\circle*{.1}}
\put(3.5,3.5){\circle{.2}}

\put(3.4,0.4){$F$}
\put(3.5,1){\circle*{.1}}
\put(3.5,1){\circle{.2}}

\put(0.2,3.3){$F'$}
\put(1,3.5){\circle*{.1}}
\put(1,3.5){\circle{.2}}

\end{picture}
\begin{center}
\captionfigure \label{5case4}
\end{center}
\end{minipage}
\hspace{2cm}
\setlength{\unitlength}{.6cm}
\begin{minipage}{4cm}
{\bf Case \ref{case5}}:\\
\begin{picture}(7,7)
\put(1,1){\framebox(1,1){}}
\put(1,2){\framebox(1,1){}}
\put(1,3){\framebox(1,1){}}
\put(1,4){\framebox(1,1){}}
\put(2,1){\framebox(1,1){}}
\put(2,2){\framebox(1,1){}}
\put(2,3){\framebox(1,1){}}
\put(2,4){\framebox(1,1){}}
\put(2,1){\framebox(1,1){}}
\put(2,2){\framebox(1,1){}}
\put(2,3){\framebox(1,1){}}
\put(2,4){\framebox(1,1){}}
\put(3,1){\framebox(1,1){}}
\put(3,2){\framebox(1,1){}}
\put(3,3){\framebox(1,1){}}
\put(3,4){\framebox(1,1){}}
\put(4,1){\framebox(1,1){}}
\put(4,2){\framebox(1,1){}}
\put(4,3){\framebox(1,1){}}
\put(4,4){\framebox(1,1){}}

\multiput(0,2.03)(.5,0){12}{\line(1,0){.3}}
\multiput(0,5.03)(.5,0){12}{\line(1,0){.3}}
\put(5.4,2.1){$h_1$}
\put(5.4,5.1){$h_2$}
\multiput(1.97,0)(0,.5){12}{\line(0,1){.3}}
\multiput(4.97,0)(0,.5){12}{\line(0,1){.3}}
\put(4.1,5.4){$v_1$}
\put(2.1,5.4){$v_2$}

\put(2.1,2.1){$A$}
\put(4.1,2.1){$B$}
\put(4.1,4.1){$C$}
\put(1.6,4.1){$D$}

\put(2,2){\circle*{.1}}
\put(2,2){\circle{.2}}
\put(4,2){\circle*{.1}}
\put(4,2){\circle{.2}}
\put(4,4){\circle*{.1}}
\put(4,4){\circle{.2}}
\put(2,4){\circle*{.1}}
\put(2,4){\circle{.2}}

\put(3.1,2.1){$P$} \put(3,2){\circle*{.25}}
\put(4.2,3.1){$Q$} \put(4,3){\circle*{.25}}
\put(2.5,4.2){$R$}\put(3,4){\circle*{.25}}
\put(1.6,3.1){$S$}\put(2,3){\circle*{.25}}

\put(1,2){\rule[-0.06\unitlength]{0.9\unitlength}{0.12\unitlength}} 
\put(2.1,2){\rule[-0.06\unitlength]{0.8\unitlength}{0.12\unitlength}} 
\put(3.1,2){\rule[-0.06\unitlength]{0.8\unitlength}{0.12\unitlength}} 
\put(4.1,2){\rule[-0.06\unitlength]{.9\unitlength}{0.12\unitlength}}

\put(3.96,2.1){\rule{0.12\unitlength}{.9\unitlength}} 
\put(3.96,4.1){\rule{0.12\unitlength}{.9\unitlength}} 
\put(3.96,1){\rule{0.12\unitlength}{.9\unitlength}} 
\put(3.1,4){\rule[-0.06\unitlength]{0.8\unitlength}{0.12\unitlength}} 
\put(1.96,2.1){\rule{0.12\unitlength}{.8\unitlength}} 
\put(1.96,3.1){\rule{0.12\unitlength}{1.9\unitlength}} 
\put(1.96,1.1){\rule{0.12\unitlength}{.9\unitlength}} 

\put(2.4,1.65){$e_1$}
\put(3.2,1.65){$e_2$}
\put(4.3,1.65){$e_3$}
\put(4.1,2.7){$e_5$}
\put(4.1,4.65){$e_6$}
\put(3.3,4.25){$e_8$}
\put(1.3,4.65){$e_{10}$}
\put(1.3,2.5){$e_{11}$}

\put(0.65,1){$0$}
\put(1,1){\circle*{.1}}
\put(1,1){\circle{.2}}

\put(3.1,3.1){$M$}
\put(3,3){\circle*{.1}}
\put(3,3){\circle{.2}}

\put(2.9,0.4){$F$}
\put(3,1){\circle*{.1}}
\put(3,1){\circle{.2}}

\put(0.4,2.9){$F'$}
\put(1,3){\circle*{.1}}
\put(1,3){\circle{.2}}

\end{picture}
\begin{center}
\captionfigure \label{6case5}
\end{center}
\end{minipage}

\end{center}
Using the same arguments one obtains:\\[3mm]
{\bf Case \ref{case2}}:\\
Removed segments: $e_1$,$e_2$, $e_3$, $e_5$, $e_6$, $e_7$, $e_{10}$, $e_{11}$.
\begin{eqnarray*}
\begin{array}{c@{\;\;,\quad}c@{\;\;,\quad}c@{\;\;,\quad}c}
l(e_1) = 0 & l(e_2) = 1 & l(e_3) = 1 & l(e_5) = 1,\\
l(e_6) = 0 & l(e_7) = 1 & l(e_{10}) = 1 & l(e_{11}) = 0.
\end{array}
\end{eqnarray*}
Thus, the leaves are switched at $e_2$, $e_3$, $e_5$, $e_7$ and $e_{10}$ .\\[3mm]
{\bf Case \ref{case3}}:\\
Removed segments: $e_1$, $e_2$, $e_5$, $e_{10}$, $e_{11}$.
\begin{eqnarray*}
\begin{array}{c@{\;\;,\quad}c@{\;\;,\quad}c@{\;\;,\quad}c@{\;\;,\quad}c}
l(e_1) = 1 & l(e_2) = 1 & l(e_5) = 1 & l(e_{10}) = 1 & l(e_{11}) = 0.
\end{array}
\end{eqnarray*}
Thus, the leaves are switched at $e_1$, $e_2$, $e_5$ and $e_{10}$.\\[3mm]
{\bf Case \ref{case4}}:\\
Removed segments: $e_1$, $e_2$, $e_4$, $e_5$.
\begin{eqnarray*}
\begin{array}{c@{\;\;,\quad}c@{\;\;,\quad}c@{\;\;,\quad}c}
l(e_1)=1 & l(e_2) = 0 &  l(e_4) = 1 & l(e_5) = 0.
\end{array}
\end{eqnarray*}
Thus, the leaves are switched at $e_1$ and $e_4$.\\[3mm]
{\bf Case \ref{case5}}:\\
Removed segments: $e_1$, $e_2$, $e_3$, $e_5$, $e_6$, $e_8$, $e_{10}$, $e_{11}$.
\begin{eqnarray*}
\begin{array}{c@{\;\;,\quad}c@{\;\;,\quad}c@{\;\;,\quad}c}
l(e_1) = 0 & l(e_2) = 1 & l(e_3) = 1 & l(e_5) =  1,\\
l(e_6) = 1 & l(e_8) = 1 & l(e_{10}) = 0 & l(e_{11}) = 1.
\end{array}
\end{eqnarray*}
Thus, the leaves are switched at $e_2$, $e_3$, $e_5$, $e_6$, $e_8$ and $e_{11}$.\\[3mm]
By the same arguments as in the first case the origamis defined by these glueings fulfill the conditions in the proposition and are uniquely determined.
\end{proof}

\end{document}